\pdfoutput=1
\documentclass[9pt,twoside,epsf]{article}
\usepackage{stmaryrd}
\usepackage{amsmath,latexsym,amssymb,amsfonts,amsbsy}
\usepackage{bm}
\usepackage{graphicx,epstopdf}
\usepackage[caption=false]{subfig}
\usepackage{cases}
\usepackage{multirow}
\usepackage{booktabs}

\newfont{\bb}{msbm10}
\def\Bbb#1{\mbox{\bb #1}}

\def\T{\top}

\def\diag{{\rm diag}}

\def\rank{{\rm rank}}

\usepackage{colortbl}
\usepackage{ragged2e}
\usepackage[margin=2em,labelsep=space,skip=0.5em,font=normalfont]{caption}
\usepackage{algorithm}  
\usepackage{algpseudocode} 
\floatname{algorithm}{\color{black} Algorithm}
\algrenewcommand{\algorithmiccomment}[1]{\quad{\color{red}\%\ #1}}
\numberwithin{algorithm}{section}
\makeatletter
\newenvironment{breakalgo}[2]{%
\captionsetup{margin=0pt,justification=RaggedRight,singlelinecheck=false}%
\par\noindent%
\medskip%
\rule{\linewidth}{0.8pt}%
\vspace{-0.5\baselineskip}%
\noindent\captionof{algorithm}{#1}\label{#2}%
\vspace{-0.7\baselineskip}%
\noindent\rule{\linewidth}{.4pt}%
\vspace{-0.3\baselineskip}%
}{%
\vspace{-.75\baselineskip}%
\rule{\linewidth}{.4pt}%
\medskip%
}

\makeatother

\newtheorem{method}{Method}[section]

\newtheorem{remark}{Remark}[section]
\newtheorem{theorem}{Theorem}[section]

\newtheorem{lemma}{Lemma}[section]

\newcommand{\reals}{\makebox{{\Bbb R}}}
\newcommand{\ceals}{\makebox{{\Bbb C}}}

\newcommand{\Z}{\makebox{{\Bbb Z}}}

\makeatletter 

\@addtoreset{equation}{section}
\makeatother 

\begin{document}
\cleardoublepage
\pagestyle{myheadings}
\bibliographystyle{plain}

\title{A fast normal splitting preconditioner for attractive coupled nonlinear Schrödinger equations with fractional Laplacian\thanks{This work was funded by the National Natural Science Foundation (No. 11101213 and No. 12071215), China. Corresponding author: Xi Yang.}}

\author{Yan Cheng\footnotemark[2] and Xi Yang\thanks{College of Mathematics, Nanjing University of Aeronautics and Astronautics, Nanjing 211106, China. Key Laboratory of Mathematical Modelling and High Performance Computing of Air Vehicles (NUAA), MIIT, Nanjing 211106, China(chengyan@nuaa.edu.cn, yangxi@nuaa.edu.cn).}}

\maketitle

\begin{abstract}
A linearly implicit conservative difference scheme is applied to discretize the attractive coupled nonlinear Schrödinger equations with fractional Laplacian. Complex symmetric linear systems can be obtained, and the system matrices are indefinite and Toeplitz-plus-diagonal. Neither efficient preconditioned iteration method nor fast direct method is available to deal with these systems. In this paper, we propose a novel matrix splitting iteration method based on a normal splitting of an equivalent real block form of the complex linear systems. This new iteration method converges unconditionally, and the quasi-optimal iteration parameter is deducted. The corresponding new  preconditioner is obtained naturally, which can be constructed easily and implemented efficiently by fast Fourier transform. Theoretical analysis indicates that the eigenvalues of the preconditioned system matrix are tightly clustered. Numerical experiments show that the new preconditioner can significantly accelerate the convergence rate of the Krylov subspace iteration methods. Specifically, the convergence behavior of the related preconditioned GMRES iteration method is spacial mesh-size-independent, and almost fractional order insensitive. Moreover, the linearly implicit conservative difference scheme in conjunction with the preconditioned GMRES iteration method conserves the discrete mass and energy in terms of a given precision.
	
	\bigskip
	
	
	\bigskip
	
	
\end{abstract}

\section{Introduction}\label{intro}
One of the most popular methods in quantum mechanics \cite{FeynmanLectures,FeynmanHSbook} is the Feynman path integral over the Brownian-like paths, which leads to the standard Schrödinger equations. By taking the path integral over the L\'{e}vy-like paths instead of the Brownian-like paths, a system of fractional Schrödinger equations (FSE) \cite{Laskin01} can be developed. Instead of the Laplacian in the standard Schrödinger equations \cite{Laskin02}, FSE includes a space fractional derivative of order $\alpha$ $(1<\alpha<2)$. When $\alpha=2$, FSE is reduced to the standard case. Besides, FSE has important physical applications in quantum mechanics, semiconductor, and other fields \cite{appli1}. Due to the non-local nature of the fractional differential operator, the exact solution of FSE is difficult to obtain. Fortunately, a large amount of numerical methods have been developed to study the nature of FSE, such as finite element methods \cite{LM2018JCP,LM2017NUMA}, spectral methods \cite{DSW2016CMA,WY2019ANM}, collocation methods \cite{Amore2010JMP,Bhrawy2017ANM}, and finite difference methods \cite{WDL2013JCP,WDL2014JCP,WPD2015JCP,ZRP2019SCM,ZX2014SISC}, etc.

In this paper, the following one dimensional (1D) space fractional coupled nonlinear Schrödinger (CNLS) equations are considered
\begin{equation} \label{equ1}
	\begin{cases}
		{\imath}u_{t}-\gamma(-\Delta)^{\frac{\alpha}{2}}u+\rho(|u|^2+\beta|v|^2)u=0,\\
		{\imath}v_{t}-\gamma(-\Delta)^{\frac{\alpha}{2}}v+\rho(|v|^2+\beta|u|^2)v=0,\
	\end{cases}
	 x \in\reals, 0 < t \le \mbox{T}
\end{equation}
with the initial conditions
	$$	u(x,0)=u_{0}(x),\quad v(x,0)=v_{0}(x), \quad x \in \reals , $$
where $\imath=\sqrt{-1}$ , the parameters $\gamma>0$, $\beta\ge 0$, $\rho$ are real constants, the fractional Laplacian  $(-\Delta)^\frac{\alpha}{2} = -\frac{\partial^\alpha}{\partial |x|^\alpha}$ is defined by the 1D Riesz space fractional derivative \cite{Laskin2002,MRiesz1949} with $1<\alpha < 2$. When the parameter $\beta=0$, the decoupled nonlinear Schrödinger (DNLS) equations \cite{BGuo2008} can be obtained. As to the parameter $\rho$, which determines the three different cases of  the nonlinear term of FSE. When $\rho=0$, the nonlinear terms are disappeared, and the system (\ref{equ1}) describes the free particles \cite{Laskin2002,Luchko2013}; when $\rho<0$, the nonlinear terms in the system (\ref{equ1}) represent the repulsive interaction of particles \cite{BWZ2012arXiv,Carr2000PRA,JS1999CPAM}; when $\rho>0$, the nonlinear terms in the system (\ref{equ1}) represent the attractive interaction of particles \cite{BWZ2012arXiv,BWZ2003SINA,Saito2001PRL}. In this paper, we only take $\rho>0$ into account, in which case the coefficient matrix of the discretized linear system is  complex symmetric indefinite. The case of $\rho<0$, leading to complex symmetric negative definite linear systems, is discussed in \cite{HSS-like09,Pihss08,RanWW2017,W-PMHSS11,ZhangYangDNTB}.

In theory, the solution of FSE keeps mass and  energy, e.g., the system (\ref{equ1}). Specifically, the solutions $u(x, t), v(x, t)$ of (\ref{equ1}) satisfy the conservation of mass \cite{WDL2014JCP}, i.e.,
$$
\begin{gathered}
	\|u(\cdot, t)\|_{L_2}^2=\int_\mathbb{R}|u(x, t)|^2 \text{d} x=\|u(\cdot, 0)\|_{L_2}^2, \\
	\|v(\cdot, t)\|_{L_2}^2=\int_\mathbb{R}|v(x, t)|^2 \text{d} x=\|v(\cdot, 0)\|_{L_2}^2 .
\end{gathered}
$$
and the conservation of energy, i.e., by defining the energy functional
$$
E(t)=\frac{\gamma}{2} \int_\mathbb{R}\left(\bar{u}(-\Delta)^{\frac{\alpha}{2}} u+\bar{v}(-\Delta)^{\frac{\alpha}{2}} v\right) \text{d} x-\frac{\rho}{4} \int_\mathbb{R} \left(|u|^4+|v|^4+2 \beta|u|^2|v|^2\right) \text{d} x,
$$
it satisfies $E(t)=E(0)$.

Due to the above conservations, the discretization scheme for the system (\ref{equ1}) should conserve mass and energy in discrete level. In this paper, we truncate the infinite space interval $\reals$ into a bounded interval $\mbox{a} \le x \le \mbox{b}$, and adopt the Dirichlet boundary conditions
\begin{align}\nonumber
		u(\mbox{a},t)=u(\mbox{b},t)=0, \quad v(\mbox{a},t)=v(\mbox{b},t)=0, \qquad  0 \le t \le \mbox{T}.
\end{align}
Then, we use a linearly implicit conservative difference (LICD) scheme proposed in \cite{WDL2014JCP} to discretize the space fractional CNLS equations (\ref{equ1}). At each time level, complex symmetric linear systems with Toeplitz-plus-diagonal structure $(D-T+\imath I)\textbf{u}=\textbf{b}$ need to be solved, where $D$ is diagonal, and $T$ is Toeplitz and symmetric positive definite. When $\rho>0$, $D$ is positive semi-definite, and $D-T$ is indefinite. Then, $D-T+\imath I$ is complex symmetric indefinite.

There are many numerical methods to deal with complex symmetric linear systems. The first class of methods are the alternating direction implicit iteration type methods \cite{HSS10,BaiGP2004,MHSS18,BaiBC2011,BaiBCW2013}. Bai et al. \cite{HSS10,BaiGP2004} proposed  the Hermitian and skew-Hermitian splitting (HSS) iteration methods to deal with general non-Hermitian positive definite linear systems. Furthermore, Bai et al.  developed the modified HSS (MHSS \cite{MHSS18}) iteration methods and the preconditioned MHSS (PMHSS \cite{BaiBC2011,BaiBCW2013}) iteration methods by considering the symmetric structure of the system matrix. The second class of methods are the C-to-R iteration methods \cite{Axelsson2000}, which require high-quality preconditioners for the Schur-complement to ensure the efficiency and robustness. The construction of the Schur-complement is difficult, since the Schur-complement matrix is dense of the form $S+\widetilde{S}^{-1}$, where $S, \widetilde{S} \in \reals^{M\times M}$ are dense matrices with Toeplitz-plus-diagonal structure. Unfortunately, the HSS type methods and the C-to-R iterations mentioned above can not be applied to handle the complex symmetric indefinite matrix $D-T+\imath I$ directly, since they are designed for non-Hermitian positive definite matrix.

For linear systems with Toeplitz-plus-diagonal structure, no fast direct solvers are developed. Hence, a better option is to adopt the preconditioned Krylov subspace iteration methods \cite{GolubBook,sad22}. Implementing matrix-vector multiplication and solving generalized residual (GR) linear system are always needed at each iteration of the preconditioned Krylov subspace iteration methods. If there are fast implementations for matrix-vector multiplication, and high-quality preconditioners can be constructed and implemented efficiently, the Krylov subspace iterations will perform really well. Many efficient preconditioners are developed to deal with the Hermitian positive definite Toeplitz-plus-diagonal systems, such as a class of banded preconditioners for Hermitian Toeplitz-plus-band systems proposed by Chan and Ng in \cite{ChanNg1993}, a class of approximate inverse circulant-plus-diagonal (AICD) preconditioners proposed by Ng and Pan in \cite{NgPan2010}, the diagonal and Toeplitz splitting (DTS) preconditioners and the diagonal and circulant splitting (DCS) preconditioners proposed by Bai et al. in \cite{BaiLuPan2017} (1D case) and \cite{BaiLu2020} (higher dimensional case), etc. Nevertheless, no preconditioners above can be directly applied to solve the complex linear systems $(D-T+\imath I)\textbf{u}=\textbf{b}$, since the performance can be poor when these preconditioners are applied to complex symmetric indefinite Toeplitz-plus-diagonal systems rather than Hermitian positive definite ones.

In order to develop effective and efficient iteration methods, we consider a real block two-by-two equivalent form of the linear system $(D-T+\imath I)\textbf{u}=\textbf{b}$, and we split the block system matrix into a normal block matrix and an anti-symmetric block matrix, and construct a class of normal and anti-symmetric splitting (NASS) iteration methods, which inherits the philosophy for designing the alternating direction implicit (ADI) iteration methods \cite{ADI1955,ADI1962}. At each step of the NASS iteration method, two linear subsystems with coefficient matrices $ \bigl[ \begin{smallmatrix} (\omega+1)I & T \\	-T & (\omega+1)I \end{smallmatrix} \bigr]$ and $ \bigl[\begin{smallmatrix} \omega I & -D \\ D & \omega I \end{smallmatrix} \bigr]$ need to be solved. The former can be solved iteratively by preconditioned GMRES (PGMRES) method with a block preconditioner of the form $ \bigl[ \begin{smallmatrix} (\omega+1)I & C \\	-C & (\omega+1)I \end{smallmatrix} \bigr]$ with a circulant approximation $C$, and each iteration of PGMRES can be accomplished in $\mathcal{O}(M\log M)$ flops if a fast algorithm is applied, i.e., fast Fourier transform (FFT) \cite{GolubBook}. The latter can be solved directly by a sparse Gaussian elimination (GE) in $\mathcal{O}(M)$ flops. In theory, the NASS iteration method converges unconditionally to the unique solution. An upper bound of the contraction factor of the NASS iteration method is derived, and it only depends on the spectra of the symmetric positive definite Toeplitz matrix $T\in\reals^{M\times M}$. The optimal value of the iteration parameter minimizing the upper bound of the contraction factor of the NASS iteration method is determined by the lower and the upper bounds of the eigenvalues of $T$.

A class of matrix splitting preconditioners, called the NASS preconditioners, can be obtained naturally from the NASS iteration method. In order to reduce the computational complexity of the NASS preconditioners, circulant approximations are used to replace the Toeplitz blocks in the NASS preconditioners, which results in a class of more economical circulant improved normal and anti-symmetric (CNAS) preconditioners. The implementation of the CNAS preconditioner requires the resolutions of two linear subsystems. The subsystem matrices are of the forms $ \bigl[ \begin{smallmatrix} (\omega+1)I & C \\	-C & (\omega+1)I \end{smallmatrix} \bigr]$ and $ \bigl[\begin{smallmatrix} \omega I & -D \\ D & \omega I \end{smallmatrix} \bigr]$. The former can be solved in $\mathcal{O}(M\log M)$ flops based on FFT, and the latter can be solved directly with sparse GE in $\mathcal{O}(M)$ flops. Theoretical analysis shows that the eigenvalues of the NASS preconditioned system matrix are clustered around 1. In addition, since it can be proved that the CNAS preconditioned system matrix is a small norm perturbation and a low rank correction of the NASS preconditioned system matrix, we can expect that most of the eigenvalues of the CNAS preconditioned system matrix are clustered around those of the NASS preconditioned system matrix. Then most of the eigenvalues of the CNAS preconditioned system matrix may also be clustered around 1. Numerical experiments show that the CNAS preconditioner improves the computational efficiency of Krylov subspace iteration methods (such as GMRES) significantly, and the corresponding PGMRES method behaves independently of the space mesh size and insensitively of the fractional order. Besides, it is observed that the discrete mass and energy remain conserved by the LICD scheme in conjunction with the CNAS preconditioned GMRES method in terms of a given precision.

The organization of the paper is as follows. In Section \ref{discre}, the LICD scheme is applied to discretize the space fractional CNLS equations, which leads to the indefinite complex symmetric linear system $(D-T+\imath I)\textbf{u}=\textbf{b}$. In Section \ref{iter}, the NASS iteration method is given, and the corresponding asymptotic convergence theory is established. In Section \ref{precon}, the new preconditioners are constructed and analyzed. In Section \ref{impcom}, the implementations and  the computational complexities of the new preconditioners are presented and discussed. In Section \ref{exp}, the numerical results are reported. Finally, in Section \ref{conclu}, we give some concluding remarks, and discuss the future work.

\section{The discretization and the linear system}\label{discre}
The space fractional CNLS equations (\ref{equ1}) are discretized by the LICD scheme \cite{WDL2014JCP}. Let $N$ and $M$ be given positive integers, we denote by $\tau=T/N$ and $h=(\text{b}-\text{a})/(M+1)$ the sizes of time steps and spatial grids, respectively. Then, $t_n=n\tau$ for $n=0,1,...,N$ and $x_j=\text{a}+jh$ for $j=0,1,...,M+1$ represent the temporal levels and the nodes for the spatial partitions.

Let $u^{n}_{j}\approx u(x_{j},t_{n})$ and  $v^{n}_{j}\approx v(x_{j},t_{n})$ be the grid function approximations. The 1D fractional Laplacian $(-\Delta)^\frac{\alpha}{2}$ can be discretized by the fractional centered difference \cite{MDuman2012,riesz21} in the bounded spatial interval $[\mbox{a},\mbox{b}]$ as
\begin{gather*}
	(-\Delta)^\frac{\alpha}{2}u(x_{j},t)
	=\frac{1}{h^{\alpha}}\Delta_h^\alpha u_{j}+\mathcal{O}(h^2),
\end{gather*}
where $\Delta_h^\alpha u_{j} =\sum^M_{k=1}c_{j-k}u_{k}$ with coefficients
\begin{gather} \nonumber
	c_{k}=(-1)^{k}\Gamma(\alpha+1)/[\Gamma(\alpha/2-k+1) \Gamma(\alpha/2+k+1)], \forall\ k\in\Z,
\end{gather}
where $\Gamma(\cdot)$ is the gamma function. In addition, the coefficients $c_{k}$ satisfy the following properties \cite{MDuman2012}
\begin{gather} \nonumber
	c_{0}\ge 0,\
	c_{k}=c_{-k}\le0\ (\forall\ k\ge 1),\ \mbox{and}\
	\sum^{+\infty}_{k=-\infty,k\ne0}|c_{k}|=c_{0}.
\end{gather}

The LICD scheme \cite{WDL2014JCP} for the system (\ref{equ1}) is as follows
\begin{equation}  \label{discretizedCNLS}
	\begin{cases}
		\imath\frac{u^{n+1}_j-u^{n-1}_j}{2\tau}-\frac{\gamma}{h^\alpha} \Delta_h^\alpha\hat{u}^n_j+\rho(|u^n_j|^2+\beta|v^n_j|^2) \hat{u}^n_j=0, \\
		\imath\frac{v^{n+1}_j-v^{n-1}_j}{2\tau}-\frac{\gamma}{h^\alpha} \Delta_h^\alpha\hat{v}^n_j+\rho(|v^n_j|^2+\beta|u^n_j|^2) \hat{v}^n_j=0,
	\end{cases}
\end{equation}
where
$\hat{u}^n_j=(u^{n+1}_j+u^{n-1}_j)/2$, $\hat{v}^n_j=(v^{n+1}_j+v^{n-1}_j)/2$, for $j=1,2,\ldots,M, n=1,2,\ldots,N-1$.
The initial and boundary conditions are $u^0_j=u_0(x_j)$, $v^0_j=v_0(x_j)$, $u^{n}_0=u^{n}_{M+1}=0$, and $v^{n}_0=v^{n}_{M+1}=0$.
Equivalently, the first line of (\ref{discretizedCNLS}) leads to a complex linear system as
\begin{align} \label{discretizedCNLSMaxForm}
	A^{n+1}\textbf{u}^{n+1} &= \textbf{b}^{n+1},\ \forall\ n\ge 1,
\end{align}
where $A^{n+1}=D^{n+1}-T+\imath I\in\ceals^{M\times M}$ , $D^{n+1}=\diag\{d^{n+1}_1,d^{n+1}_2,\ldots,d^{n+1}_M\}\in\reals^{M\times M}$ is diagonal with $d^{n+1}_j=\rho\tau(|u^n_j|^2+\beta |v^n_j|^2)$, $I\in\reals^{M\times M}$ is the identity matrix, and $T\in\reals^{M\times M}$ is Toeplitz of the form
\begin{align}\label{equ5}
	T &=\mu
	\begin{bmatrix}
		c_0 & c_{-1} & \ldots & c_{2-M} & c_{1-M} \\
		c_1 & c_0 & \ddots & \ddots & c_{2-M} \\
		\vdots & \ddots & \ddots & \ddots & \vdots \\
		c_{M-2} & \ddots & \ddots & c_0 & c_{-1} \\
		c_{M-1} & c_{M-2} & \ldots & c_1 & c_0 \\
	\end{bmatrix}
\end{align}
with $\mu=\frac{\gamma\tau}{h^\alpha}$. Clearly, the second line of (\ref{discretizedCNLS}) also admits a coupled linear system the same as   (\ref{discretizedCNLSMaxForm}).

Due to the facts that $\rho>0$, $\gamma>0$, $\beta \ge 0$ and the properties of the coefficients $c_k$, $D^{n+1}$ is a nonnegative diagonal matrix, $T$ is strictly diagonally dominant and symmetric. Hence, $D^{n+1}-T$ is symmetric indefinite, and $A^{n+1}$ is complex symmetric and indefinite.

It is proved in \cite{WDL2014JCP} that the discretized coupled linear systems (\ref{discretizedCNLS}) conserve the discrete mass and energy, i.e.,
$$
Q_1^n=Q_1^0, \quad Q_2^n=Q_2^0, \quad E^n=E^0, \quad 1 \leq n \leq N,
$$
where $Q_1^n=(\|u^{n+1}\|^2+\|u^n\|^2)/2$ and $Q_2^n=(\|v^{n+1}\|^2+\|v^n\|^2)/2$ with $\|u\|^2=\langle u, u\rangle$ and $\langle u, v\rangle=h \sum_{k=1}^{M-1} u_j \bar{v}_j$, and
$$
\begin{aligned}
	E^n= & \frac{\gamma h}{4 h^\alpha} \sum_{j=1}^{M-1}\left(\bar{u}_j^{n+1} \Delta_h^\alpha u_j^{n+1}+\bar{u}_j^n \Delta_h^\alpha u_j^n+\bar{v}_j^{n+1} \Delta_h^\alpha v_j^{n+1}+\bar{v}_j^n \Delta_h^\alpha v_j^n\right) \\
	 &-\frac{\rho h}{4} \sum_{j=1}^{M-1}\left(\left(\left|u_j^n\right|^2\left|u_j^{n+1}\right|^2+\left|v_j^n\right|^2\left|v_j^{n+1}\right|^2\right)+\beta\left(\left|u_j^n\right|^2\left|v_j^{n+1}\right|^2+\left|v_j^n\right|^2\left|u_j^{n+1}\right|^2\right)\right) .
\end{aligned}\label{energyformula}
$$

\section{The NASS iteration method}\label{iter}
We consider the construction of iteration method for the complex linear systems of the form (\ref{discretizedCNLSMaxForm}). By ignoring the superscripts, it reads
\begin{align} \label{equ3}
	A \textbf{u} &= \textbf{b},
\end{align}
where $A=D-T+\imath I\in\ceals^{M\times M}$ is complex symmetric and indefinite, with $D \in\reals^{M \times M}$ a positive semi-definite diagonal matrix and $T \in \reals^{M \times M}$ a symmetric positive definite Toeplitz matrix. The complex unknown vector is of the form $\textbf{u}=y+\imath z\in\ceals^{M}$, and the complex vector $\textbf{b}=p+\imath q\in\ceals^{M}$ is the right-hand-side, where $y$, $z$, $p$, $q\in\reals^{M}$ are real vectors. Then, the complex linear system (\ref{equ3}) can be reformed to the following real block linear system
\begin{align}\label{equ11}
	\widehat{\mathcal{R}}\widehat{x} &\equiv
	\begin{bmatrix}
		D-T & -I\\
		I & D-T\\
	\end{bmatrix}
	\begin{bmatrix}
		y \\
		z\\
	\end{bmatrix}
	=
	\begin{bmatrix}
		p \\
		q \\
	\end{bmatrix}
	\equiv \widehat{f},
\end{align}
where $\widehat{\mathcal{R}}\in \reals^{2M\times 2M}$ is non-symmetric and indefinite.

By defining the block permutations $\mathcal{P} = \bigl[\begin{smallmatrix}-I & 0 \\ 0 & I\end{smallmatrix}\bigr]$, $\mathcal{Q} = \bigl[\begin{smallmatrix} 0 & I \\ I & 0 \end{smallmatrix} \bigr]$,
and let
\begin{align} \nonumber \mathcal{R}=\mathcal{P}\widehat{\mathcal{R}}\mathcal{Q},\ x=\mathcal{Q}\widehat{x},\ f=\mathcal{P}\widehat{f},\end{align}
the real block linear system (\ref{equ11}) is equivalent to a real non-symmetric positive definite block linear system as follows
\begin{align}\label{positiveBlockForm}
	\mathcal{R}x &\equiv
	\begin{bmatrix}
		I & T-D\\
		D-T & I\\
	\end{bmatrix}
	\begin{bmatrix}
		z \\
		y\\
	\end{bmatrix}
	=
	\begin{bmatrix}
		-p \\
		q \\
	\end{bmatrix}
	\equiv f.
\end{align}

The system matrix $\mathcal{R}\in\reals^{2M\times 2M}$ admits the following normal and anti-symmetric (NAS) splitting
\begin{align} \label{NASSsplitting}
	\mathcal{R} &= \mathcal{T} + \mathcal{D},
\end{align}
where $\mathcal{T}= \bigl[ \begin{smallmatrix} I & T \\
	-T & I \end{smallmatrix} \bigr]$ is normal, $\mathcal{D}= \bigl[ \begin{smallmatrix} 0 & -D \\ D & 0 \end{smallmatrix} \bigr]$ is anti-symmetric. Based on the NAS splitting  (\ref{NASSsplitting}) and the spirit of the alternating direction implicit (ADI) iteration \cite{ADI1955,ADI1962}, we can construct the NASS iteration method as follows.

\begin{method}[The NASS iteration method]
	\label{NASSiteration}
	Let $x^{(0)}\in\reals^{2M} $ be an arbitrary initial guess. For $k=0,1,2,\ldots$ until the sequence of iterates $\{x^{(k)}\}_{k\ge 0}$ converges, compute the next iterate $x^{(k+1)}$ according to the following procedure:
	\begin{align}\label{equ12}
		\left\{\begin{aligned}
			(\omega I+\mathcal{T})x^{(k+\frac{1}{2})} & = (\omega I-\mathcal{D})x^{(k)}+f, \\
			(\omega I+\mathcal{D})x^{(k+1)} & = (\omega I-\mathcal{T})x^{(k+\frac{1}{2})}+f,
		\end{aligned}\right.
	\end{align}
	where $\omega>0$ is an arbitrary iteration parameter.
	\end{method}
	
	The iteration scheme (\ref{equ12}) can be reformulated to a fixed point iteration as
	$x^{(k+1)} = \mathcal{L}_{\omega}x^{(k)}+\mathcal{F}^{-1}_{\omega}f$,
	where the iteration matrix reads
	\begin{align} \label{NASSiterMatrix}
		\mathcal{L}_{\omega} &= \mathcal{F}^{-1}_{\omega}\mathcal{G}_{\omega},
	\end{align}
	where
	\begin{align} \label{defpreconditionger1}
		\mathcal{F}_{\omega} = \frac{1}{2\omega}
		(\omega I+\mathcal{T})(\omega I+\mathcal{D}), \end{align}
	and
	\begin{align} \nonumber
		\mathcal{G}_{\omega} = \frac{1}{2\omega}
		(\omega I-\mathcal{T})(\omega I-\mathcal{D}).
	\end{align}
	In fact, the matrices $\mathcal{F}_{\omega}$ and $\mathcal{G}_{\omega}$ form a splitting of $\mathcal{R}$, i.e.,
	\begin{align} \label{matrixSplittingR}
		\mathcal{R} &= \mathcal{F}_{\omega}-\mathcal{G}_{\omega}.
	\end{align}
	
	In the following theorem, the convergence property of the NASS iteration method \ref{NASSiteration} is provided.
	
	\begin{theorem} \label{NASSconvergenceThm}
		Let $\mathcal{R}\in\reals^{2M\times 2M}$ be the system matrix in (\ref{positiveBlockForm}). Let $\mathcal{T}\in\reals^{2M\times 2M}$ and $\mathcal{D}\in\reals^{2M\times 2M}$ be the block matrices forming the splitting of  $\mathcal{R}$ in (\ref{NASSsplitting}). Let $\omega$ be a positive constant. Then, the NASS iteration method (\ref{equ12}) is well defined, and the iteration matrix $\mathcal{L}_{\omega}$ is given by (\ref{NASSiterMatrix}). The spectral radius $\rho(\mathcal{L}_{\omega})$ is bounded by
		\begin{align} \label{NASSupperBound}
			\sigma(\omega) &= \max_{\lambda_i\in\lambda(T)} \sqrt{\frac {(\omega-1)^2+\lambda_i^2}{(\omega+1)^2+\lambda_i^2}},
		\end{align}
		where $\lambda(T)$ is the spectral set of $T$. Therefore, it reads
		\begin{align} \label{NASSrateBoundedByOne}
			\rho(\mathcal{L}_{\omega}) \le \sigma(\omega) < 1,\ \forall\ \omega > 0,
		\end{align}
		i.e., the NASS iteration converges to the unique solution of the block linear system (\ref{positiveBlockForm}).
	\end{theorem}
	
	{\em Proof.}
	Since $\mathcal{T}$ is positive definite, $\mathcal{D}$ is anti-symmetric, and $\omega>0$, it follows that $\omega I+\mathcal{T}$ and $\omega I+\mathcal{D}$ are invertible and positive definite. Thus, the NASS iteration (\ref{equ12}) is well defined. Note that the iteration matrix $\mathcal{L}_{\omega}$ satisfies
	$\mathcal{L}_{\omega}
	= (\omega I+ \mathcal{D})^{-1}\mathcal{U}_{\omega} \mathcal{V}_{\omega}(\omega I +\mathcal{D})$,
	where $\mathcal{U}_{\omega} = (\omega I+\mathcal{T})^{-1}(\omega I-\mathcal{T})$ and $\mathcal{V}_{\omega} = (\omega I-\mathcal{D})(\omega I+\mathcal{D})^{-1}$, the following relation can be established
	\begin{align} \nonumber
		\rho(\mathcal{L}_{\omega}) &= \rho(\mathcal{U}_{\omega} \mathcal{V}_{\omega}) \\ \label{NASSrateBoundedBy2norm}
		&\le \|\mathcal{U}_{\omega}\|_2 \|\mathcal{V}_{\omega}\|_2.
	\end{align}
	
	Now, we focus on estimating $\|\mathcal{U}_{\omega}\|_2$ and $\|\mathcal{V}_{\omega}\|_2$.
	
	Firstly, it is easy to verify that $\mathcal{T}^{\T}\mathcal{T}=\mathcal{T}\mathcal{T}^{\T}$, then we have
	\begin{align} \nonumber
		\mathcal{U}_{\omega}^{\T}\mathcal{U}_{\omega} &= (\omega I -\mathcal{T})^{\T}(\omega I +\mathcal{T})^{-\T}(\omega I +\mathcal{T})^{-1}(\omega I -\mathcal{T}) \\ \nonumber
		&= (\omega I +\mathcal{T})^{-\T}(\omega I +\mathcal{T})^{-1}(\omega I -\mathcal{T})^{\T}(\omega I -\mathcal{T}) \\ \nonumber
		&= [(\omega I +\mathcal{T})(\omega I +\mathcal{T})^{\T}]^{-1}(\omega I -\mathcal{T})^{\T}(\omega I -\mathcal{T})\\
		&= \begin{bmatrix}
			(\omega +1)^2I+T^2& 0 \\
			0 & (\omega +1)^2I+T^2 \
		\end{bmatrix}^{-1}
		\begin{bmatrix}
			(\omega -1)^2I+T^2 & 0 \\
			0 & (\omega -1)^2I+T^2 \
		\end{bmatrix}.
	\end{align}
	Due to the fact that $T$ is symmetric positive definite, it holds that
	\begin{align}\label{norm-2-u}
		\|\mathcal{U}_{\omega}\|_2 &=  \max_{\lambda_i\in\lambda(T)} \sqrt{\frac {(\omega-1)^2+\lambda_i^2}{(\omega+1)^2+\lambda_i^2}} < 1.
	\end{align}
	
	Secondly, we consider the estimate of $\|\mathcal{V}_{\omega}\|_2$. Similarly, since  $\mathcal{D}^{\T}\mathcal{D}=\mathcal{D}\mathcal{D}^{\T}$, it follows that
	\begin{align*}
		\mathcal{V}_{\omega}^{\T}\mathcal{V}_{\omega} &= (\omega I+\mathcal{D}^{\T})^{-1} (\omega I-\mathcal{D}^{\T}) (\omega I-\mathcal{D}) (\omega I+\mathcal{D})^{-1} \\
		&= (\omega I-\mathcal{D}^{\T}) (\omega I-\mathcal{D}) [(\omega I+\mathcal{D}) (\omega I+\mathcal{D}^{\T})]^{-1} \\
		&=
		\begin{bmatrix}
			\omega^2 I+D^2 & 0 \\
			0 & \omega^2 I+D^2 \
		\end{bmatrix}
		\begin{bmatrix}
			\omega^2 I+D^2 & 0 \\
			0 & \omega^2 I+D^2 \
		\end{bmatrix}^{-1}\\
		&= I.
	\end{align*}
	Then, it reads
	\begin{align} \label{NASSupBoundFactorV}
		\left \|\mathcal{V}_{\omega}\right \|_2 &=   1.
	\end{align}
	
	Based on the facts (\ref{NASSrateBoundedBy2norm})-(\ref{NASSupBoundFactorV}), the estimate (\ref{NASSrateBoundedByOne}) can be achieved.
	$\hfill\square$

	\begin{remark}
		To further understand the theoretical result stated in the above theorem, we provide the following remarks.
		\begin{enumerate}
			\item The convergence rate of the NASS iteration is bounded by $\sigma(\omega)$,  only depending on the spectra of the Toeplitz matrix $T$.
			\item If we introduce a weighted vector norm $\|x\|_{\mathcal{D}}=\|(\omega I+\mathcal{D}) x\|_2$ for $x\in\reals^{2M}$ and the corresponding induced weighted matrix norm $\|\mathcal{X}\|_{\mathcal{D}}=\|(\omega I+\mathcal{D}) \mathcal{X} (\omega I+\mathcal{D})^{-1}\|_2$ for $\mathcal{X}\in\reals^{2M\times 2M}$,  it holds that
			\begin{align*}
				\|\mathcal{L}_{\omega}\|_{\mathcal{D}} = \|\mathcal{U}_{\omega}\mathcal{V}_{\omega}\|_2 < \sigma(\omega),\ \forall\ \omega>0.
			\end{align*}
			Moreover, based on the fixed point iteration form of the NASS iteration, it can be verified that the sequence of the iterates $\{x^{(k)}\}_{k\ge 0}$ satisfies
			\begin{align*}
				\|x^{(k+1)} - x^{\star}\|_{\mathcal{D}} \le \sigma(\omega) \|x^{(k)} - x^{\star}\|_{\mathcal{D}},\ \forall\ k\ge 0.
			\end{align*}
			Hence, $\sigma(\omega)$ also serves as an upper bound of the contraction factor of the NASS iteration in the sense of the $\|\cdot\|_{\mathcal{D}}$-norm. It is noted that, when $TD=DT$, it holds $\mathcal{U}_{\omega} \mathcal{V}_{\omega} =\mathcal{V}_{\omega} \mathcal{U}_{\omega}$, leading to a fact $\rho(\mathcal{L}_{\omega}) =\| \mathcal{L}_{\omega} \|_{\mathcal{D}}= \sigma({\omega})$. Then, these quantities can be minimized at the same optimal value of $\omega$.
			\item The upper bound $\sigma(\omega)$ can be further bounded by $\widehat{\sigma}(\omega)$, i.e.,
			\begin{align} \nonumber
				\sigma(\omega) & \le  \max_{\lambda_{\min}\le\lambda\le\lambda_{\max}} \sqrt{\frac {(\omega-1)^2+\lambda^2}{(\omega+1)^2+\lambda^2}}
				\equiv\widehat{\sigma}(\omega),
			\end{align}
			where $\lambda_{\min}>0$ and $\lambda_{\max}>0$ are the minimum and the maximum eigenvalues of $T$. Since the function $g(\omega;\lambda)=\sqrt{[(\omega-1)^2+\lambda^2]/[(\omega+1)^2+\lambda^2]}$ is monotonically increasing for $\lambda>0$, it reads
			\begin{align} \nonumber
				\widehat{\sigma}(\omega) &=  \sqrt{\frac{(\omega-1)^2+\lambda_{\max}^2}{(\omega+1)^2+\lambda_{\max}^2}}.
			\end{align}
			Let
		$$
				\frac{\rm{d}(\widehat{\sigma}(\omega))}{\rm{d} \omega}=\frac{2\left(\omega^2-1-\lambda_{\max }^2\right)}{\left[(\omega-1)^2+\lambda_{\max }^2\right]^{\frac{1}{2}}\left[(\omega+1)^2+\lambda_{\max }^2\right]^{\frac{3}{2}}}=0,$$
then $\widehat{\sigma}(\omega)$ is minimized at $\omega^{\star}=\sqrt{\lambda_{\max}^2+1}$. If the parameter $\omega=\omega^{\star}$ is applied, the convergence rate of the NASS iteration for solving the block linear system (\ref{positiveBlockForm}) is bounded as follows
			\begin{align*}
				\sigma(\omega^{\star}) \le \widehat{\sigma}(\omega^{\star}) =\frac{\lambda_{\max}}{1+\sqrt{\lambda_{\max}^2+1}}.
			\end{align*}
			
		\end{enumerate}
	\end{remark}

\section{Preconditioning}\label{precon}
The NASS iteration method naturally induces the preconditioning matrix $\mathcal{F}_{\omega}$ defined in (\ref{defpreconditionger1}) for the system matrix $\mathcal{R}$ in (\ref{positiveBlockForm}).
The corresponding preconditioned linear system reads
\begin{align} \label{preconPositiveBlockForm}
	\mathcal{F}_{\omega}^{-1} \mathcal{R}x &= \mathcal{F}_{\omega}^{-1} f.
\end{align}
The preconditioner $\mathcal{F}_{\omega}$ is called the NASS preconditioner. Obviously, $\mathcal{F}_{\omega}$ is the product of a scalar $1/(2\omega)$, and the block normal matrices $\omega I+\mathcal{T}$ and $\omega I+\mathcal{D}$.

In actual implementation, the main workload of the NASS preconditioning is to solve a sequence of GR linear systems
\begin{align} \label{generalizedResEqns}
	\mathcal{F}_{\omega}z^{(k)} = r^{(k)},\ \forall\ k\ge 0,
\end{align}
where $r^{(k)}$ is the current residual vector, and $z^{(k)}$ is the GR vector. The resolution of (\ref{generalizedResEqns}) consists of solving the linear subsystems with coefficient matrices $\omega I+\mathcal{T}$ and $\omega I+\mathcal{D}$.
In order to reduce the cost of the preconditioning process, circulant approximations are considered to replace the Toeplitz matrix $T$ in $\mathcal{F}_{\omega}$ \cite{CRH1989SISSC, CT1988SISSC, CRH1996SIREV}. In theory, we only consider the Strang's circulant approximation $C$ for $T$  \cite{CRH1989SISSC}.
For simplicity of the analysis, let $M$ be even in the sequel, it reads $C=\mu\ [\varsigma_{ij}] \in\reals^{M\times M}$ with $\varsigma_{ij}=\varsigma_{i-j}$, $\varsigma_{k} = \varsigma_{-k}$ for $0\le k< M$, $\varsigma_{\frac{M}{2}} =0$, and
\begin{align*}
	\varsigma_{k} &= \begin{cases}
		c_k, & \mbox{if } 0\le k < \frac{M}{2}, \\
		c_{M-k}, & \mbox{if } \frac{M}{2} < k < M.
	\end{cases}
\end{align*}
i.e.,
\begin{align}\label{cform}
	C=\mu\left[\begin{array}{cccccccc}
		c_0 & c_{-1} & \cdots & c_{-\left(\frac{M}{2}-1\right)} & 0 & c_{-\left(\frac{M}{2}-1\right)} & \cdots & c_{-1} \\
		c_{-1} & c_0 & \ddots & \ddots & c_{-\left(\frac{M}{2}-1\right)} & 0 & \ddots & \vdots \\
		\vdots & \ddots & \ddots & \ddots & \ddots & \ddots & \ddots & c_{-\left(\frac{M}{2}-1\right)} \\
		c_{-\left(\frac{M}{2}-1\right)} & \ddots & \ddots & \ddots & \ddots & \ddots & \ddots & 0 \\
		0 & \ddots & \ddots & \ddots & \ddots & \ddots & \ddots & c_{-\left(\frac{M}{2}-1\right)} \\
		c_{-\left(\frac{M}{2}-1\right)} & 0 & \ddots & \ddots & \ddots & \ddots & \ddots & \vdots \\
		\vdots & \ddots & \ddots & \ddots & \ddots & \ddots & \ddots & c_{-1} \\
		c_{-1} & c_{-2} & \cdots & 0 & c_{-\left(\frac{M}{2}-1\right)} & \ddots & \cdots & c_0
	\end{array}\right],
\end{align}
which is real symmetric. Then,  an improved efficient circulant based preconditioner $\widetilde{\mathcal{F}}_{\omega}$ is built, i.e.,
\begin{align*}
	\widetilde{\mathcal{F}}_{\omega} &= \frac{1}{2\omega}
	(\omega I+\mathcal{C})(\omega I+\mathcal{D}),
\end{align*}
where $\mathcal{C}= \bigl[ \begin{smallmatrix} I & C \\
	-C & I \end{smallmatrix} \bigr]$. Obviously, $\mathcal{C} $ is   normal, and $\mathcal{D}$ is anti-symmetric. Thus, the preconditioner $\widetilde{\mathcal{F}}_{\omega}$ is called the circulant improved normal and anti-symmetric (CNAS) preconditioner.

Now, we study the clustering property of the eigenvalues of the preconditioned system matrix  $\widetilde{\mathcal{F}}^{-1}_\omega \mathcal{R}$. Due to the fact that
\begin{align} \label{CNASpreconSysMax}
	\widetilde{\mathcal{F}}^{-1}_\omega \mathcal{R} &= \underbrace{\widetilde{\mathcal{F}}^{-1}_\omega \mathcal{F}_{\omega}}\ \underbrace{\mathcal{F}_{\omega}^{-1} \mathcal{R}},
\end{align}
we discuss the properties of $\mathcal{F}_{\omega}^{-1} \mathcal{R}$ and $\widetilde{\mathcal{F}}^{-1}_\omega \mathcal{F}_{\omega}$ respectively.

 Firstly, we consider the clustering property of $\mathcal{F}_{\omega}^{-1} \mathcal{R}$.
\begin{theorem} \label{propertyNASSpreconCoeffiMax}
	Let $\mathcal{R}\in\reals^{2M\times 2M}$ be the system matrix in (\ref{positiveBlockForm}), $\mathcal{T}\in\reals^{2M\times 2M}$ and $\mathcal{D}\in\reals^{2M\times 2M}$ be the block matrices forming the splitting of  $\mathcal{R}$ in (\ref{NASSsplitting}), $\omega$ be a positive constant, and $\sigma(\omega)$ be defined by (\ref{NASSupperBound}). When the NASS  preconditioner $\mathcal{F}_{\omega}$ is applied to the block linear system (\ref{positiveBlockForm}), the eigenvalues of the preconditioned system matrix $\mathcal{F}_{\omega}^{-1} \mathcal{R}$ in (\ref{preconPositiveBlockForm}) are located in a circle of radius $\sigma(\omega)<1$ centered at 1.
\end{theorem}

{\em Proof.}
According to (\ref{NASSiterMatrix}) and (\ref{matrixSplittingR}), it reads $\mathcal{L}_{\omega} = I - \mathcal{F}_{\omega}^{-1} \mathcal{R}$.
Thus, let $\eta$ be an eigenvalue of $\mathcal{F}_{\omega}^{-1} \mathcal{R}$, then $1-\eta$ is an eigenvalue of $\mathcal{L}_{\omega}$. In addition, Theorem \ref{NASSconvergenceThm} shows that $\rho(\mathcal{L}_{\omega})\le \sigma(\omega)$, i.e., $|1-\eta| \le \sigma(\omega)$ for all $\eta\in\lambda (\mathcal{F}_{\omega}^{-1} \mathcal{R})$, which is exactly the result of this theorem.
$\hfill\square$

\begin{theorem} \label{presystemseigen0}
	Let $\mathcal{R}\in\reals^{2M\times 2M}$ be the system matrix in (\ref{positiveBlockForm}), $\mathcal{T}\in\reals^{2M\times 2M}$ and $\mathcal{D}\in\reals^{2M\times 2M}$ be the block matrices forming the splitting of $\mathcal{R}$ in (\ref{NASSsplitting}), $\omega$ be a positive constant. Then, the eigenvalues of the preconditioned system matrix $\mathcal{F}_{\omega}^{-1} \mathcal{R}$ are clustered at $0_+$ on the right-half complex plane if $\omega$ tends to the positive infinity.
	
\end{theorem}

{\em Proof.}
Let $\lambda$ be an eigenvalue of $\mathcal{L}_{\omega}$, and $x$ be the corresponding unit eigenvector, then it reads $\mathcal{L}_{\omega}x = \lambda x$, i.e.,
\begin{align}\nonumber
	(\omega I - \mathcal{T})(\omega I - \mathcal{D})x = \lambda(\omega I + \mathcal{T})(\omega I + \mathcal{D})x.
\end{align}
Denote by $\mathcal{\widehat{T}} = \bigl[\begin{smallmatrix}0 & T \\ -T & 0\end{smallmatrix}\bigr]$, it holds that $\mathcal{T} = I + \mathcal{\widehat{T}}$, then we have
\begin{align}\nonumber
	[(\omega-1) I - \mathcal{\widehat{T}}](\omega I - \mathcal{D})x = \lambda[(\omega+1) I + \mathcal{\widehat{T}}](\omega I + \mathcal{D})x.
\end{align}
By multiplying $x^*$ from left on both sides of the above equation, it follows that
\begin{align}\nonumber
\lambda = \frac{\omega(\omega-1)-\omega x^*\mathcal{\widehat{T}}x-(\omega-1)x^*\mathcal{D}x+x^*\mathcal{\widehat{T}D}x}{\omega(\omega+1)+\omega x^*\mathcal{\widehat{T}}x+(\omega+1)x^*\mathcal{D}x+x^*\mathcal{\widehat{T}D}x}.
\end{align}
Since $ x^*\mathcal{\widehat{T}}x$ and $x^*\mathcal{D}x$ are pure imaginary numbers due to the facts that $\mathcal{\widehat{T}}$ and $\mathcal{D}$ are real and anti-symmetric, we denote by $ x^*\mathcal{\widehat{T}}x={\gamma}_{_T}\imath$ and $x^*\mathcal{D}x=\gamma_{_D}\imath$ with ${\gamma}_{_T}, \gamma_{_D} \in \reals$. In addition, we denote by $x^*\mathcal{\widehat{T}D}x={\gamma}_{_{TD}}^{(R)}+\imath {\gamma}_{_{TD}}^{(I)}$ with ${\gamma}_{_{TD}}^{(R)}, {\gamma}_{_{TD}}^{(I)} \in \reals$. Then, it reads
\begin{align} \nonumber
	\lambda &= \frac{\omega(\omega-1)+{\gamma}_{_{TD}}^{(R)}+\imath [-\omega {\gamma}_{_T}-(\omega-1){\gamma}_{_D}+ {\gamma}_{_{TD}}^{(I)}]}{\omega(\omega+1)+{\gamma}_{_{TD}}^{(R)}+\imath [\omega {\gamma}_{_T}+(\omega+1){\gamma}_{_D}+ {\gamma}_{_{TD}}^{(I)}]}.
\end{align}
Since $\mathcal{F}_{\omega}^{-1} \mathcal{R} = I-\mathcal{L}_\omega$, then 1-$\lambda$ is an eigenvalue of $\mathcal{F}_{\omega}^{-1} \mathcal{R}$, and it reads
\begin{align}\nonumber
	1-\lambda &= \frac{2}{\omega}\frac{1+\imath({\gamma}_{_T}+{\gamma}_{_D})+\mathcal{O}(1/\omega)}{[1+1/\omega+{\gamma}_{_{TD}}^{(R)}/{\omega^2}]^2+[{\gamma}_{_T}/\omega+(\omega+1){\gamma}_{_D}/{\omega^2}+{\gamma}_{_{TD}}^{(I)}/{\omega^2}]^2}.
\end{align}
Obviously, the real part of 1-$\lambda$ tends to $0_+$, and the imaginary part of 1-$\lambda$ tends to 0, when the parameter $\omega$ tends to the positive infinity.

$\hfill\square$

Secondly, we discuss the property of $\widetilde{\mathcal{F}}^{-1}_\omega \mathcal{F}_{\omega}$. For convenience, we define the following constants
$$\theta=\frac{\left(1-\frac{1+\alpha}{5+\alpha / 2}\right)^{5+\frac{\alpha}{2}} e^{1+\alpha} \Gamma(\alpha+1) \sin \left(\frac{\pi \alpha}{2}\right)}{\pi \alpha} \quad \text { and } \quad \theta_0=\frac{\sqrt{2} e^{13 / 12} \Gamma(\alpha+1) \sin \left(\frac{\pi \alpha}{2}\right)}{\pi \alpha}.$$
Then, two lemmas are introduced as follows.
\begin{lemma}[\cite{ZhangYangDNTB}]\label{lemma1} Let $c_j=(-1)^j \Gamma(\alpha+1) /[\Gamma(\alpha / 2-j+1) \Gamma(\alpha / 2+j+1)], k_0 \geq 3$, and $1<\alpha < 2$, then
	$$
\frac{\theta}{(k_0+1/2)^{\alpha}}<\sum_{j=k_0+1}^{\infty}\left|c_j\right|<\frac{\theta_0}{(k_0-1)^{\alpha}}.
	$$
\end{lemma}

Based on Lemma \ref{lemma1}, the bounds for the eigenvalues of the Toeplitz matrix $T$ and its circulant approximation $C$ are presented in the following lemma.

\begin{lemma}[\cite{ZhangYangDNTB}]\label{lemma2}
	Let $T$ be the Toeplitz matrix in (\ref{equ5}), $C$ be the Strang's circulant approximation in (\ref{cform}), and $M$ be even, then it holds that
	$$
	\frac{2 \gamma \tau \theta}{(\mathrm{b}-\mathrm{a})^\alpha}<\lambda_T<\frac{2 \gamma \tau}{h^\alpha}\left[\frac{\Gamma(\alpha+1)}{\Gamma(\alpha / 2+1)^2}-\frac{\theta h^\alpha}{(\mathrm{b}-\mathrm{a})^\alpha}\right], \quad M \geq 4,
	$$
	and
	$$
	\frac{2^{\alpha+1} \gamma \tau \theta}{(\mathrm{b}-\mathrm{a})^\alpha}<\lambda_C<\frac{2 \gamma \tau}{h^\alpha}\left[\frac{\Gamma(\alpha+1)}{\Gamma(\alpha / 2+1)^2}-\frac{2^\alpha \theta h^\alpha}{(\mathrm{b}-\mathrm{a})^\alpha}\right], \quad M \geq 8,
	$$
	where $\lambda_T$ and $\lambda_C$ are the eigenvalues of $T$ and $C$.
\end{lemma}

Based on Lemmas \ref{lemma1} and \ref{lemma2}, the property of ${\mathcal{\widetilde{F}}}^{-1}_\omega \mathcal{F}_{\omega}$ is summarized in the following theorem.
\begin{theorem} \label{propertyTildeFinvF}
	Let $1<\alpha < 2$, and $M\ge8$ be even. Let $\epsilon>0$ be a small constant, and $k_0=\left\lceil\left(\frac{\mu \theta_0}{\epsilon}\right)^{\frac{1}{\alpha}}\right\rceil+1$, where $\lceil$.$\rceil $ represents rounding a real number to the positive infinity. Then, there exist two matrices $\mathcal{\widetilde{E}} \in \mathbb{R}^{2 M \times 2 M}$ and $\mathcal{\widetilde{F}} \in \mathbb{R}^{2 M \times 2 M}$, satisfying $\operatorname{rank}\left(\mathcal{\widetilde{E}}\right)=4 k_0$, $$
	\left\|\mathcal{\widetilde{E}}\right\|_2 <\frac{(\omega^2 + \nu^2)^{\frac{1}{2}}M^{\frac{1}{2}} \mu}{\omega\sqrt{(\omega+1)^2+[\frac{2^{\alpha+1} \gamma \tau \theta}{(\rm{b}-\rm{a})^\alpha}]^2}}\left[\frac{c_0}{2}-\frac{\theta}{\left(M-\frac{1}{2}\right)^\alpha}\right] \quad \text { and } \quad\left\|\mathcal{\widetilde{F}}\right\|_2 < \frac{(\omega^2 + \nu^2)^{\frac{1}{2}}M^{\frac{1}{2}} \epsilon}{\omega\sqrt{(\omega+1)^2+[\frac{2^{\alpha+1} \gamma \tau \theta}{(\rm{b}-\rm{a})^\alpha}]^2}},
	$$
where $\nu = \max_{\mu_i \in \lambda(D)} |\mu_i|$, such that
\begin{align} \label{tildeFinvF}
	\widetilde{\mathcal{F}}^{-1}_\omega \mathcal{F}_{\omega} &= I + \mathcal{\widetilde{E}}+\mathcal{\widetilde{F}}.
\end{align}
\end{theorem}

{\em Proof.}
Since $M$ is even, it reads
\begin{gather*}
	T-C=\mu
	\begin{bmatrix}
		0 & \widehat{F}_{12} & \widehat{E}_{13} \\
		\widehat{F}^{\T}_{12} & 0 & 0 \\
		\widehat{E}^{\T}_{13} & 0 & 0 \\
	\end{bmatrix},
\end{gather*}
where $\widehat{F}_{12}\in \reals^{\frac{M}{2}\times(\frac{M}{2}-k_0)}$ and $\widehat{E}_{13}\in \reals^{\frac{M}{2}\times k_0}$ are of the forms
\begin{align*}
	\widehat{F}_{12}&=
	\begin{bmatrix}
		c_{\frac{M}{2}} & c_{\frac{M}{2}+1}-c_{\frac{M}{2}-1} & \cdots & c_{M-(k_0+1)}-c_{k_0+1} \\
		0 & \ddots & \ddots & \vdots \\
		\vdots & \ddots & \ddots & c_{\frac{M}{2}+1}-c_{\frac{M}{2}-1} \\
		\vdots & \ddots & \ddots & c_{\frac{M}{2}} \\
		0 & \cdots & \cdots & 0 \\
		\vdots &   &   & \vdots \\
		0 & \cdots & \cdots & 0
	\end{bmatrix}, \\
	\widehat{E}_{13}&=
	\begin{bmatrix}
		c_{M-k_0}-c_{k_0} & \cdots & \cdots & c_{M-1}-c_1 \\
		\vdots & \ddots & \ddots & \vdots \\
		c_{\frac{M}{2}+1}-c_{\frac{M}{2}-1} & \ddots & \ddots & \vdots \\
		c_{\frac{M}{2}} & \ddots & \ddots & c_{M-k_0}-c_{k_0} \\
		0 & \ddots & \ddots & \vdots \\
		\vdots & \ddots & \ddots & c_{\frac{M}{2}+1}-c_{\frac{M}{2}-1} \\
		0 & \cdots & 0 & c_{\frac{M}{2}}
	\end{bmatrix}.
\end{align*}

Therefore, it holds that \begin{align}\label{t-c}
	T-C=\widehat{E}+\widehat{F},\end{align}
where
\begin{gather*}
	\widehat{E}=\mu
	\begin{pmatrix}
		0 & 0 & \widehat{E}_{13} \\
		0 & 0 & 0\\
		\widehat{E}^{\T}_{13} & 0 & 0 \\
	\end{pmatrix}\ \mbox{and}\
	\widehat{F}=\mu
	\begin{pmatrix}
		0 & \widehat{F}_{12} & 0 \\
		\widehat{F}^{\T}_{12} & 0 & 0 \\
		0 & 0 & 0\\
	\end{pmatrix}.
\end{gather*}
Obviously, it reads $\rank(\widehat{E})=2k_0$. Moreover, due to Lemmas \ref{lemma1} and \ref{lemma2}, and the structure of  $\widehat{E}$ and $\widehat{F}$, the following estimates can be obtained, i.e.,
$$
\begin{aligned}
	\|\widehat{E}\|_{\infty} & =\mu \max \left\{\left\|\widehat{E}_{13}\right\|_{\infty},\left\|\widehat{E}_{13}^{\top}\right\|_{\infty}\right\}=\mu\left\|\widehat{E}^{\T}_{13}\right\|_{\infty} \\
	& \leq \mu \sum_{j=1}^{M-1}\left|c_j\right|=\mu\left(\frac{c_0}{2}-\sum_{j=M}^{\infty}\left|c_j\right|\right) \\
	& <\mu\left[\frac{c_0}{2}-\frac{\theta}{\left(M-\frac{1}{2}\right)^\alpha}\right],
\end{aligned}
$$

$$
\begin{aligned}
	\|\widehat{F}\|_{\infty} & =\mu \max \left\{\left\|\widehat{F}_{12}\right\|_{\infty},\left\|\widehat{F}_{12}^{\top}\right\|_{\infty}\right\}=\mu\left\|\widehat{F}_{12}\right\|_{\infty} \\
	& \leq \mu \sum_{j=k_0+1}^{M-k_0-1} c_j<\mu \sum_{j=k_0+1}^{\infty}\left|c_j\right| \\
	& <\frac{\mu \theta_0}{\left(k_0-1\right)^\alpha}<\epsilon.
\end{aligned}
$$

It can be verified that  $\widetilde{\mathcal{F}}^{-1}_\omega \mathcal{F}_{\omega} - I$ and $(\omega I + \mathcal{C})^{-1} (\mathcal{T} - \mathcal{C})$ are similar, i.e.,
\begin{align} \nonumber
	\widetilde{\mathcal{F}}^{-1}_\omega \mathcal{F}_{\omega}-I &=
	(\omega I+\mathcal{D})^{-1}\left[(\omega I+\mathcal{C})^{-1}(\omega I+\mathcal{T})-I\right](\omega I+\mathcal{D}) \\
	&= (\omega I+\mathcal{D})^{-1}(\omega I+\mathcal{C})^{-1}(\mathcal{T}-\mathcal{C})(\omega I+\mathcal{D}), \nonumber
\end{align}
with $(\omega I + \mathcal{C})^{-1} (\mathcal{T} - \mathcal{C}) =
\bigl[ \begin{smallmatrix}	(\omega+1)I & C \\
	-C & (\omega+1)I \end{smallmatrix} \bigr]^{-1}
	\bigl[ \begin{smallmatrix}	0 & T-C \\
		C-T & 0 \end{smallmatrix} \bigr]$.
Together with (\ref{t-c}), it reads
\begin{align} \nonumber
	(\omega I + \mathcal{C})^{-1} (\mathcal{T} - \mathcal{C}) &= \mathcal{\widehat{E}} + \mathcal{\widehat{F}},
\end{align}
where
\begin{align} \nonumber
	\mathcal{\widehat{E}}= \begin{bmatrix}
		(\omega+1)I & C \\
		-C & (\omega+1)I
	\end{bmatrix}^{-1}\begin{bmatrix}
		0 & \widehat{E}\\
		-\widehat{E} & 0
	\end{bmatrix} \intertext{and}          \nonumber
	\mathcal{\widehat{F}}=\begin{bmatrix}
		(\omega+1)I & C \\
		-C & (\omega+1)I
	\end{bmatrix}^{-1}\begin{bmatrix}
		0 & \widehat{F}\\
		-\widehat{F} & 0
	\end{bmatrix}.
\end{align}
Therefore, it holds that
\begin{align} \nonumber
	\widetilde{\mathcal{F}}^{-1}_\omega \mathcal{F}_{\omega}-I & =
	\mathcal{\widetilde{E}}+ \mathcal{\widetilde{F}},
\end{align}
where
$	\mathcal{\widetilde{E}}= (\omega I+\mathcal{D})^{-1} 	\mathcal{\widehat{E}} (\omega I+\mathcal{D})$ and $	\mathcal{\widetilde{F}} = (\omega I+\mathcal{D})^{-1} 	\mathcal{\widehat{F}}(\omega I+\mathcal{D})$.

In addition, we have
\begin{align} \nonumber
	\left \|\mathcal{\widetilde{E}}\right \|_2 &\le  \left\|\left[\begin{array}{cc}\omega I & -D \\ D & \omega I\end{array}\right]^{-1}\right\|_2  \left\|\left[\begin{array}{cc}\omega I & -D \\ D & \omega I\end{array}\right]\right\|_2   \left\|\left[\begin{array}{cc}(\omega+1) I & C \\ -C & (\omega+1)I\end{array}\right]^{-1}\right\|_2\|\widehat{E}\|_2 \\ \nonumber
	&= \frac{\max_{\mu_i\in\lambda(D)} \sqrt{\omega^2+\mu_i ^2}} {\min_{\mu_i\in\lambda(D)} \sqrt{\omega^2+\mu_i ^2}} \frac{1}{\min_{\lambda_i^{(c)}  \in\lambda(C)}\sqrt{(\omega+1)^2+(\lambda_i^{(c)})^2}}\|\widehat{E}\|_2  \\ \nonumber &\leq\frac{\max_{\mu_i\in\lambda(D)} \sqrt{\omega^2+\mu_i ^2}} {\min_{\mu_i\in\lambda(D)} \sqrt{\omega^2+\mu_i ^2}} \frac{1}{\min_{\lambda_i^{(c)} \in\lambda(C)}\sqrt{(\omega+1)^2+(\lambda_i^{(c)})^2}} M^{\frac{1}{2}}\|\widehat{E}\|_{\infty} \\ \nonumber
	&<\frac{(\omega^2 + \nu^2)^{\frac{1}{2}}M^{\frac{1}{2}} \mu}{\omega\sqrt{(\omega+1)^2+[\frac{2^{\alpha+1} \gamma \tau \theta}{(\rm{b}-\rm{a})^\alpha}]^2}}\left[\frac{c_0}{2}-\frac{\theta}{\left(M-\frac{1}{2}\right)^\alpha}\right],
	\end{align}
	the first ``$=$'' is due to the facts that the eigenvalues of $  \bigl[\begin{smallmatrix}\omega I & -D \\ D & \omega I\end{smallmatrix}\bigr]$ are  $\omega \pm  \mu_i\imath$ with $\mu_i \in \lambda(D)$ and $\imath=\sqrt{-1}$, and the eigenvalues of  $\bigl[\begin{smallmatrix}(\omega+1)I  & C \\ -C & (\omega+1) I\end{smallmatrix}\bigr]$ are $\omega+1\pm\lambda_i^{(c)}\imath$ with $\lambda_i^{(c)} \in \lambda(C)$, the second ``$\le$'' is because of the relationship between $\|$·$\|_2$ and  $\|$·$\|_\infty$, the last ``$<$'' is due to Lemma \ref{lemma2}, the estimate of 	$\|\widehat{E}\|_{\infty}$ and the facts $\omega > 0$ and $0 \le \mu_i \le \nu$. Similarly, we have
	\begin{align}\nonumber
	\left \|\mathcal{\widetilde{F}}\right \|_2
	&< \frac{(\omega^2 + \nu^2)^{\frac{1}{2}}M^{\frac{1}{2}} \epsilon}{\omega\sqrt{(\omega+1)^2+[\frac{2^{\alpha+1} \gamma \tau \theta}{(\rm{b}-\rm{a})^\alpha}]^2}}.
\end{align}
$\hfill\square$

\begin{remark} \label{rmk:Preconditioning1}
	We consider the following relation stated in Theorem \ref{propertyTildeFinvF}, i.e.,
	\begin{align*}
		\widetilde{\mathcal{F}}^{-1}_\omega \mathcal{F}_{\omega} &= \underbrace{I + \mathcal{\widetilde{F}}}+\mathcal{\widetilde{E}}.
	\end{align*}
	On one hand, the matrix $I + \mathcal{\widetilde{F}}$ is a small perturbation of the identity matrix $I$. Specifically, since $\|\mathcal{\widetilde{F}}\|_2\le \epsilon$, the Bauer-Fike theorem \cite{BauerFike1960} leads to a fact that
	\begin{align*}
		|\xi-1| & \le \epsilon,\ \forall\ \xi \in \lambda(I + \mathcal{\widetilde{F}}),
	\end{align*}
	i.e., the eigenvalues of $I + \mathcal{\widetilde{F}}$ are clustered in a small disk centered at 1 with radius $\epsilon$. On the other hand, since the matrix $\mathcal{\widetilde{E}}$ has bounded $\ell_2$-norm and low rank, the matrix $\widetilde{\mathcal{F}}^{-1}_\omega \mathcal{F}_{\omega}$ can be considered as a low rank modification of $I + \mathcal{\widetilde{F}}$. Then, one can expect that the eigenvalues of $\widetilde{\mathcal{F}}^{-1}_\omega \mathcal{F}_{\omega}$ are also located in the same disk centered at 1 with radius $\epsilon$ except for a small number of outliers.
\end{remark}

Finally, the following theorem states a relationship between $\widetilde{\mathcal{F}}_{\omega}^{-1}\mathcal{R}$ and $\mathcal{F}_{\omega}^{-1}\mathcal{R}$.

\begin{theorem} \label{th442}  Let $1<\alpha<2$, $M \geq 8$ be even,
 $\epsilon$ be a small positive constant satisfying $\frac{2^\alpha \mu \theta_0}{(M-2)^\alpha}<\epsilon \leq \mu \theta_0$,  $k_0=\left\lceil\left(\frac{\mu \theta_0}{\epsilon}\right)^{\frac{1}{\alpha}}\right\rceil+1$, and $\nu = \max_{\mu_i \in \lambda(D)} |\mu_i|$. Then, there exist two matrices  $\mathcal{P}_\omega \in \mathbb{R}^{2 M \times 2 M}$ and $\mathcal{Q}_\omega \in \mathbb{R}^{2 M \times 2 M}$, satisfying $\operatorname{rank}\left(\mathcal{P}_\omega\right)=4 k_0,  \left\|\mathcal{P}_\omega\right\|_2 < \frac{2(\omega^2 + \nu^2)M^{\frac{1}{2}} \mu}{\omega^2\sqrt{(\omega+1)^2+[\frac{2^{\alpha+1} \gamma \tau \theta}{(\rm{b}-\rm{a})^\alpha}]^2}}\left[\frac{c_0}{2}-\frac{\theta}{\left(M-\frac{1}{2}\right)^\alpha}\right]$, and $\left\|\mathcal{Q}_\omega\right\|_2 < \frac{2(\omega^2 + \nu^2)M^{\frac{1}{2}} \epsilon}{\omega^2\sqrt{(\omega+1)^2+[\frac{2^{\alpha+1} \gamma \tau \theta}{(\rm{b}-\rm{a})^\alpha}]^2}}$, such that
$$
\widetilde{\mathcal{F}}_\omega^{-1} \mathcal{R}=\mathcal{F}_\omega^{-1} \mathcal{R}+\mathcal{P}_\omega+\mathcal{Q}_\omega.
$$
\end{theorem}

{\em Proof.}
It is noted that
\begin{align} \nonumber
	\mathcal{F}_\omega^{-1} \mathcal{R} & =I-\mathcal{F}_\omega^{-1} \mathcal{G}_\omega \\
	& \nonumber=\left[\begin{array}{cc}
		\omega I & -D\\
		D & \omega I
	\end{array}\right]^{-1}\left(I-\mathcal{U}_{\omega}\mathcal{V}_{\omega}\right)\left[\begin{array}{cc}
		\omega I& -D \\
		D & \omega I
	\end{array}\right],
\end{align}
where $\mathcal{U}_{\omega}$ and $\mathcal{V}_{\omega}$ are defined in the proof of Theorem \ref{NASSconvergenceThm}. According to the facts (\ref{norm-2-u}) and (\ref{NASSupBoundFactorV}), we know that $\left\|\mathcal{U}_{\omega}\mathcal{V}_{\omega}\right\|_2<1$, then $\left\|I-\mathcal{U}_{\omega}\mathcal{V}_{\omega}\right\|_2 \leq\|I\|_2+\left\|\mathcal{U}_{\omega}\mathcal{V}_{\omega}\right\|_2<2$.

Moreover, according to the proof of Theorem \ref{propertyTildeFinvF}, we have
\begin{align} \nonumber
\left\|  \left[\begin{array}{cc}
	\omega I & -D\\
	D & \omega I
\end{array}\right]^{-1}    \right\|_2\left\|  \left[\begin{array}{cc}
\omega I& -D \\
D & \omega I
\end{array}\right]    \right\|_2
&< \frac{\sqrt{\omega^2 + \nu^2}}{\omega}.
\end{align}
Then, it reads $\left\|\mathcal{F}_\omega^{-1} \mathcal{R}\right\|_2 < \frac{2\sqrt{\omega^2 + \nu^2}}{\omega}$.

Due to the relations (\ref{CNASpreconSysMax}) and (\ref{tildeFinvF}), it holds that
$$
\begin{aligned}
	\widetilde{\mathcal{F}}_\omega^{-1} \mathcal{R} & =\left(I+\mathcal{\widetilde{E}}+\mathcal{\widetilde{F}}\right) \mathcal{F}_\omega^{-1} \mathcal{R} \\
	& =\mathcal{F}_\omega^{-1} \mathcal{R}+\mathcal{P}_\omega+\mathcal{Q}_\omega
\end{aligned}
$$
with $\mathcal{P}_\omega=\mathcal{\widetilde{E}} \mathcal{F}_\omega^{-1} \mathcal{R}$ and $\mathcal{Q}_\omega=\mathcal{\widetilde{F}} \mathcal{F}_\omega^{-1} \mathcal{R}$. By Theorem \ref{propertyTildeFinvF}, it follows that $\operatorname{rank}\left(\mathcal{P}_\omega\right)=4 k_0$,
\begin{align}
	\left\|\mathcal{P}_\omega\right\|_2 & \nonumber \leq\left\|\mathcal{\widetilde{E}}\right\|_2\left\|\mathcal{F}_\omega^{-1} \mathcal{R}\right\|_2 \\
	&  \nonumber <   \frac{2(\omega^2 + \nu^2)M^{\frac{1}{2}} \mu}{\omega^2\sqrt{(\omega+1)^2+[\frac{2^{\alpha+1} \gamma \tau \theta}{(b-a)^\alpha}]^2}}\left[\frac{c_0}{2}-\frac{\theta}{\left(M-\frac{1}{2}\right)^\alpha}\right],
\intertext{and}\nonumber
	\left\|\mathcal{Q}_\omega\right\|_2 & \leq\left\|\mathcal{\widetilde{F}}\right\|_2\left\|\mathcal{F}_\omega^{-1} \mathcal{R}\right\|_2 \\
	&  \nonumber <  \frac{2(\omega^2 + \nu^2)M^{\frac{1}{2}} \epsilon}{\omega^2\sqrt{(\omega+1)^2+[\frac{2^{\alpha+1} \gamma \tau \theta}{(\rm{b}-\rm{a})^\alpha}]^2}}.
\end{align}
$\hfill\square$

\begin{remark} \label{rmk:Preconditioning2}
	Theorem \ref{th442} shows that $\mathcal{Q}_\omega$ is a small norm matrix, so the eigenvalues of $\mathcal{F}_\omega^{-1} \mathcal{R}+\mathcal{Q}_\omega$ can be considered as small perturbations of the eigenvalues of $\mathcal{F}_\omega^{-1} \mathcal{R}$. In addition, $\mathcal{P}_\omega$ has bounded $\ell_2$-norm and low rank, so we can expect that most of the eigenvalues of $	\widetilde{\mathcal{F}}_\omega^{-1} \mathcal{R} $  are distributed near the eigenvalues of  $\mathcal{F}_\omega^{-1} \mathcal{R}+\mathcal{Q}_\omega$. In summary, the eigenvalues of $\widetilde{\mathcal{F}}_\omega^{-1} \mathcal{R} $ are clustered around those of $\mathcal{F}_\omega^{-1} \mathcal{R}$ except for several outliers.
\end{remark}

\section{Implementation and complexity}\label{impcom}
In the previous sections, we proposed and analyzed new preconditioners for the block linear system (\ref{positiveBlockForm}), i.e., the NASS preconditioner $\mathcal{F}_\omega$ and the CNAS preconditioner $\widetilde{\mathcal{F}}_{\omega}$. In this section, we discuss the implementation of these preconditioners in detail. In fact, we can  ignore the scalar factor $1/(2\omega)$ of $\mathcal{F}_\omega$ and $\widetilde{\mathcal{F}}_{\omega}$   for the reason that it does not affect the property of the preconditioned system matrices. Then, the simplified preconditioners read
\begin{align} \nonumber
\mathcal{F}_{\mbox{\tiny NASS}}  &=
	\begin{bmatrix}
		\widehat{\omega}I & T \\
		-T & \widehat{\omega}I
	\end{bmatrix}
	\begin{bmatrix}
		\omega I  & -D \\
		D & \omega I
	\end{bmatrix} \\
	&\nonumber=\begin{bmatrix}
		I & 0 \\
		\widehat{L}_{21} & I
	\end{bmatrix}
	\begin{bmatrix}
		\widehat{\omega}I & T \\
		0 & \widehat{U}_{22}
	\end{bmatrix}
	\begin{bmatrix}
		I & 0 \\
		L_{21} & I
	\end{bmatrix}
	\begin{bmatrix}
		\omega I & -D \\
		0 & U_{22}
	\end{bmatrix}
	\intertext{and} \nonumber
	\mathcal{F}_{\mbox{\tiny CNAS}} &=
	\begin{bmatrix}
	\widehat{\omega}I  & C \\
		-C & \widehat{\omega}I
	\end{bmatrix}
	\begin{bmatrix}
		\omega I  & -D \\
		D & \omega I
	\end{bmatrix}\\
	&\nonumber=
	\begin{bmatrix}
		F & 0\\
		0 & F
	\end{bmatrix}^{-1}
	\begin{bmatrix}
		I & 0 \\
	\widetilde{L}_{21}& I
	\end{bmatrix}
	\begin{bmatrix}
		\widehat{\omega}I & \Lambda \\
		0 & \widetilde{U}_{22}
	\end{bmatrix}
	\begin{bmatrix}
		F & 0\\
		0 & F
	\end{bmatrix}
	\begin{bmatrix}
		I & 0 \\
		{L_{21}} & I
	\end{bmatrix}
	\begin{bmatrix}
		\omega I & -D \\
		0 & U_{22}
	\end{bmatrix},
\end{align}
where $\omega>0$ serves as the parameter of $\mathcal{F}_{\mbox{\tiny NASS}}$ and $\mathcal{F}_{\mbox{\tiny CNAS}}$, $F \in \mathbb{C}^{M \times M}$ is the discrete Fourier transform (DFT), $\Lambda=\operatorname{diag}(F c) \in \mathbb{C}^{M \times M}$ with $c \in \mathbb{C}^M$ being the first column of $C$,  $\widehat{\omega}=\omega+1$, $\widehat{L}_{21}=-T/\widehat{\omega}$,  $\widehat{U}_{22}=\widehat{\omega}I+T^2/\widehat{\omega}$, $L_{21}=D/\omega$, $U_{22}=\omega I + D^2/\omega$, $\widetilde{L}_{21}=-\Lambda/\widehat{\omega}$, $\widetilde{U}_{22}=\widehat{\omega}I+\Lambda^2/\widehat{\omega}$. The diagonal matrices $L_{21}, U_{22}, \widetilde{L}_{21}, \widetilde{U}_{22}$ can be computed in advance.

In the preconditioned Krylov subspace methods, the GR vector is computed by solving the GR equation associated with the new preconditioners at each iteration. Let $x=(x_1^{\T},x_2^{\T})^{\T}\in\reals^{2M}$ with $x_1$, $x_2\in\reals^M$, and
$r=(r_1^{\T},r_2^{\T})^{\T}\in\reals^{2M}$ with $r_1$, $r_2\in\reals^M$, the details for implementing the NASS preconditioning and the CNAS preconditioning are described in Algorithms \ref{alg-NASS}-\ref{alg-CNAS}.

The cost of Algorithm \ref{alg-CNAS} consists of $M$-vector FFT/IFFT operations, and $M$-vector operations (including vector addition, diagonal matrix-vector multiplication, and diagonal linear system solve). According to \cite{CRH2007SIAMbook}, FFT/IFFT operations can be accomplished in $\mathcal{O}(M\text{log}M)$ flops. Additionally, all the $M$-vector operations can be accomplished in $\mathcal{O}(M)$ flops. Therefore, the workload for implementing $\mathcal{F}_{\mbox{\tiny CNAS}}$ is dominated by the FFT/IFFT operations. It means that the GR vector can be computed in $\mathcal{O}(M\text{log}M)$ flops at each iteration of the preconditioned Krylov subspace methods.

The cost of Algorithm \ref{alg-NASS} consists of $M$-vector operations (around $\mathcal{O}(M)$ flops), $M\times M$ Toeplitz-matrix-vector multiplications (around $\mathcal{O}(M\text{log}M)$ flops since it can be implemented based on FFT/IFFT), and one $M\times M$ dense linear system solve with coefficient matrix $(\omega+1)I+T^2/(\omega+1)$. The above dense linear system solve can be implemented based on a direct method (around $\mathcal{O}(M^3)$), or the preconditioned conjugate gradient (PCG) method with preconditioner  $(\omega+1)I+C^2/(\omega+1)$ (around $\mathcal{O}(kM\text{log}M)$ flops, where $k$ is the iteration counts of PCG). Obviously, the workload for implementing  $\mathcal{F}_{\mbox{\tiny NASS}}$ is dominated by the cost of solving the dense linear system with coefficient matrix  $(\omega+1)I+T^2/(\omega+1)$, which is much more expensive than the workload for implementing $\mathcal{F}_{\mbox{\tiny CNAS}}$. Therefore, the new preconditioner $\mathcal{F}_{\mbox{\tiny CNAS}}$  is highly recommended.

\begin{breakalgo}{Solve the GR equation  $\mathcal{F}_{\mbox{\tiny NASS}}\,x=r$}{alg-NASS}
	\begin{algorithmic}[1]
		\State $x_1 = r_1$, $x_2 = r_2 - \widehat{L}_{21}r_1$; \% \texttt{One $M\times M$ Toeplitz-matrix-vector multiplication, two $M$-vector operations}
        \State Solve $\widehat{U}_{22}x_2 = x_2$ and $x_1 = (r_1-Tx_2)/\widehat{\omega}$; \% \texttt{One $M\times M$ dense linear system solve, one $M\times M$ Toeplitz-matrix-vector multiplication, two $M$-vector operations}
		\State $x_2 = x_2 - L_{21}x_1$; \% \texttt{Two $M$-vector operations}
		\State Solve $U_{22}x_2 = x_2$ and $x_1 = (x_1+D x_2)/\omega$. \% \texttt{Four $M$-vector operations}
	\end{algorithmic}
\end{breakalgo}

\begin{breakalgo}{Solve the GR equation  $\mathcal{F}_{\mbox{\tiny CNAS}}\,x=r$}{alg-CNAS}
\begin{algorithmic}[1]
	\State $x_j$ = FFT($r_j$), $j=1, 2$; \% \texttt{Two $M$-vector FFTs}
	\State $x_2 = x_2 - \widetilde{L}_{21}x_1$; \% \texttt{Two $M$-vector operations}
	\State Solve $\widetilde{U}_{22}x_2 = x_2$ and $x_1 = (x_1-\Lambda x_2)/\widehat{\omega}$; \% \texttt{Four $M$-vector operations}
	\State $x_j$ = IFFT($x_j$), $j = 1, 2$; \% \texttt{Two $M$-vector IFFTs}
	\State $x_2 = x_2 - L_{21}x_1$; \% \texttt{Two $M$-vector operations}
	\State Solve $U_{22}x_2 = x_2$ and $x_1 = (x_1+D x_2)/\omega$. \% \texttt{Four $M$-vector operations}
\end{algorithmic}
\end{breakalgo}

\section{Numerical experiments}\label{exp}
In this section, we present the evidence of the properties of the new preconditioners $\mathcal{F}_{\mbox{\tiny NASS}}$ and $\mathcal{F}_{\mbox{\tiny CNAS}}$, and we apply the preconditioned GMRES method with the CNAS preconditioner $\mathcal{F}_{\mbox{\tiny CNAS}}$ to solve the discretized space fractional CNLS equations in the case of attractive interaction of particles. The discretized space fractional CNLS equations are derived from the LICD scheme. The LICD scheme requires the initial value at the initial time level and a second or higher order approximate value at the first time level to start up. For instance, the second order approximate value can be obtained by a second order implicit conservative scheme \cite{WDL2013JCP}.

In all the numerical experiments, the block linear system (\ref{positiveBlockForm}) at the second time level of the discretized fractional CNLS equations is selected to be the tested linear system, and the initial guess of the preconditioned GMRES method is set to be zero vector. In addition, the (preconditioned) GMRES method is running without restart, and terminates either the $\ell_2$-norm relative residual of the tested linear system reduced below $10^{-6}$ or the iteration counts exceeding $3000$.

\subsection{The DNLS case}
Let $\beta=0$, then the system (\ref{equ1}) is decoupled. We consider the following truncated system
\begin{align}\label{equ:DNLS}
	\imath u_t-\gamma(-\Delta)^{\frac{\alpha}{2}}u + \rho\vert u \vert^2u=0, \qquad  -20\le x \le 20,\quad 0 < t \le \mbox{T},
\end{align}
subjected to the initial and boundary conditions
\begin{align}\label{equ:DNLS_InitBdry}
	u(x,0)=\text{sech}(x) \  e^{2\imath x},\ u(-20,t)=u(20,t)=0.
\end{align}
Here, we take the parameters $\gamma=1$, $\rho=2$, $1<\alpha\le 2$. The discretized space fractional DNLS equations are obtained by applying the LICD scheme to (\ref{equ:DNLS}) and (\ref{equ:DNLS_InitBdry}). It requires to solve a complex symmetric linear system of the form (\ref{equ3}) on each time level $t_n$ for $1 < n \le N$, which is equivalent to solve a block linear system (\ref{positiveBlockForm}). Specifically, the  block linear system (\ref{positiveBlockForm}) is solved by the CNAS preconditioned GMRES (CNAS-GMRES) method. In this subsection, `IT' represents the iteration counts of a tested iteration method.

\begin{figure}[htbp]
	\centering
	\includegraphics[scale=0.4]{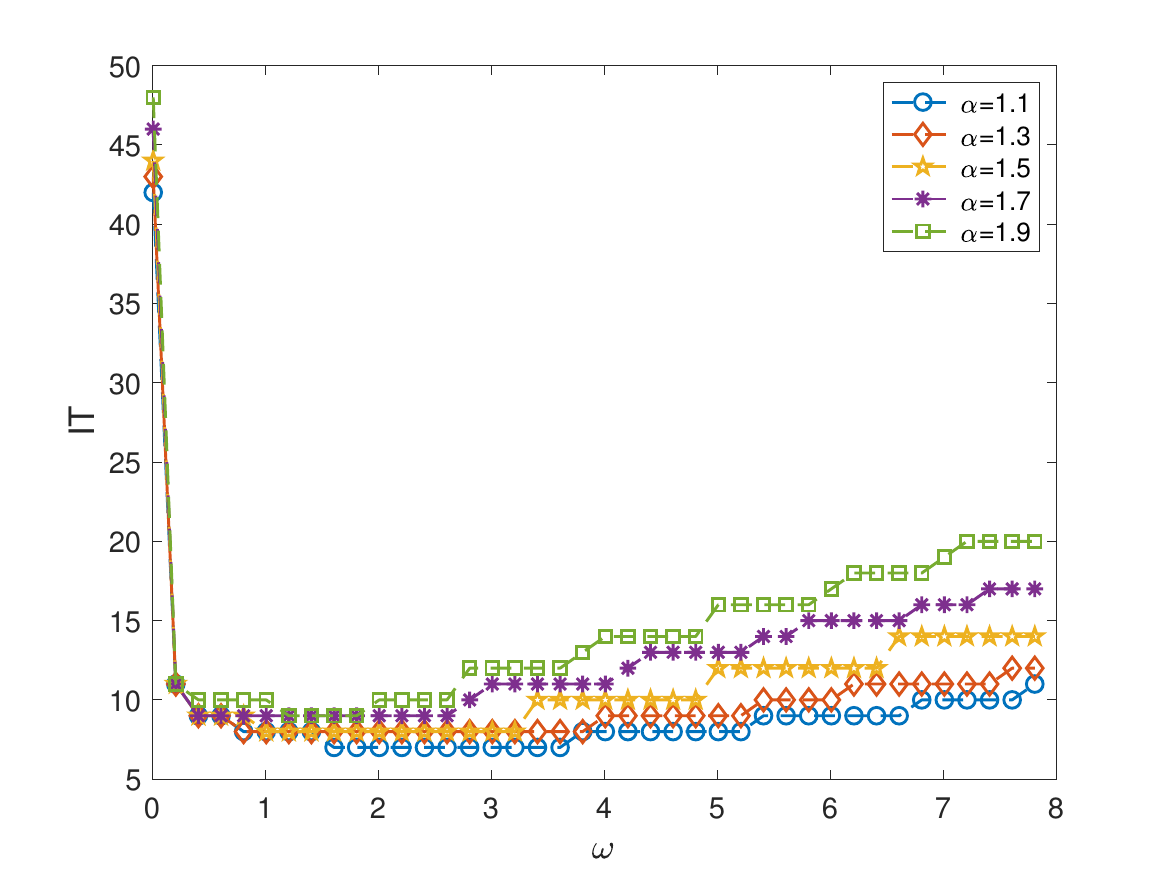}
	\caption{
		The curves of  IT versus the parameter $\omega\in (0,8]$ of CNAS-GMRES when $\alpha=1.1:0.2:1.9$, $M=6400, N=200$.
	}
	\label{fig:singal_W-it}
\end{figure}

Figure \ref{fig:singal_W-it} depicts the curves of IT versus the parameter $\omega\in (0,8]$ of CNAS-GMRES when $\alpha=1.1:0.2:1.9$, $M=6400, N=200$. From Figure \ref{fig:singal_W-it}, we can see that IT increases rapidly as $\omega$ goes to zero. However, when $\omega$ grows, IT reaches its minimum quickly and then grows very slowly, and the optimal value of $\omega$ is roughly around $[2.2,2.8]$ for all cases of $\alpha$. In summary, the convergence of CNAS-GMRES is less sensitive to the parameter $\omega$ as long as $\omega$ is not too close to zero, which makes the CNAS preconditioner $\mathcal{F}_{\mbox{\tiny CNAS}}$ easy to use. In addition, a larger $\alpha$ leads to a larger IT, which means that the system is more difficult to solve as $\alpha$ increases.

\begin{figure}[htbp]
	\centering
	\includegraphics[scale=0.4]{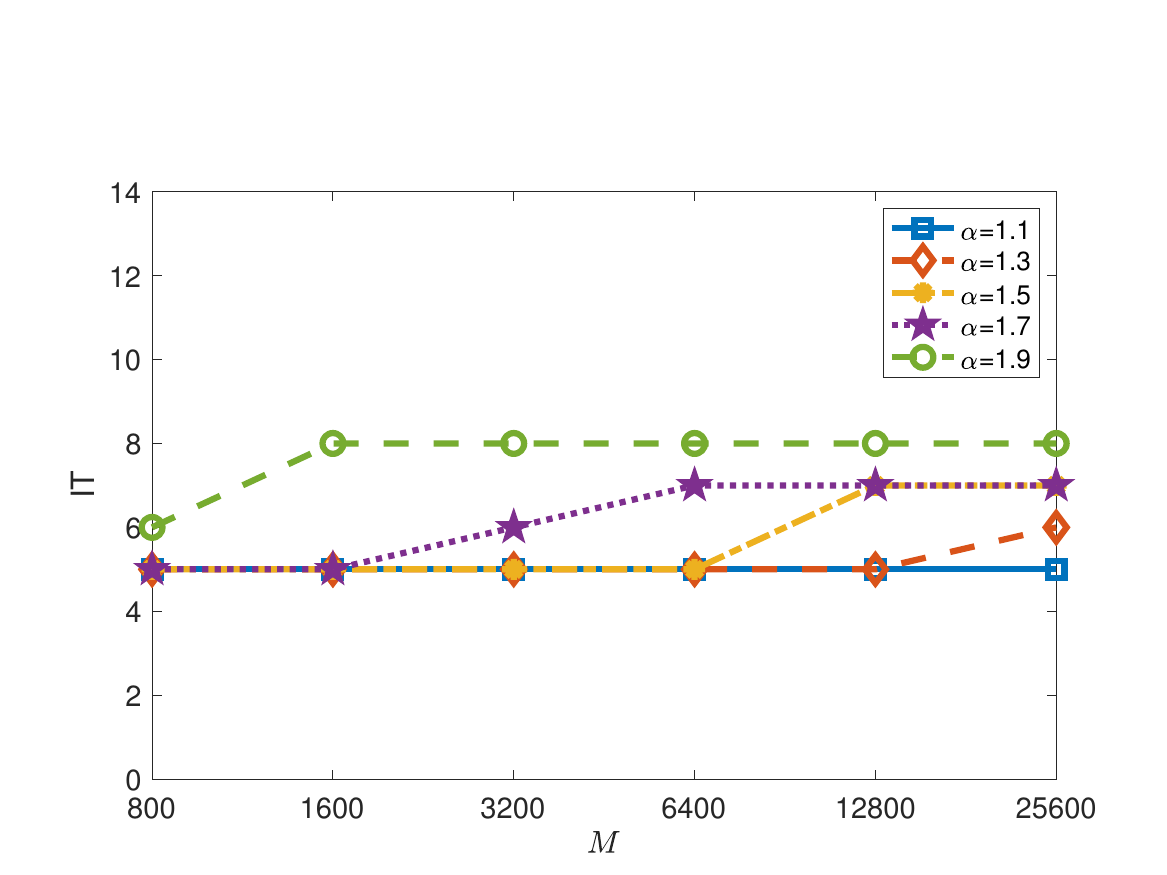}
	\caption{
		The curves of IT of CNAS-GMRES versus the space mesh size $M$ when $\alpha=1.1:0.2:1.9$, $N=200$.
	}
	\label{figitm}
\end{figure}

Figure \ref{figitm} depicts the curves of IT of CNAS-GMRES versus the space mesh size $M$ of the LICD scheme applied to the space fractional DNLS equations when $\alpha=1.1:0.2:1.9$, $N=200$. The empirical optimal value of $\omega$ is selected. It shows that IT of CNAS-GMRES remains almost the same with the increase of the number of spatial discrete points, which indicates that the convergence property of CNAS-GMRES is independent of the space mesh size. In addition, a larger fractional order $\alpha$ leads to a higher position of the curve in the plot.

Table \ref{tab:TNS_different_circulant_alp=1.9} lists IT of CNAS-GMRES in conjunction with different circulant matrices, and the empirical optimal value of the parameter $\omega$ of CNAS-GMRES when $\alpha=1.9$, $M=6400, N=200$. The search interval of the parameter $\omega$ is $(0,4]$. Only the results of several representative circulant matrices are listed \cite{CRH2007SIAMbook}, including  T. Chan's circulant matrix, Strang's circulant matrix, R. Chan's circulant matrix, circulant matrices constructed from some famous kernels (e.g., Modified Dirichlet kernel, von Hann kernel, Hamming kernel), and superoptimal circulant matrix. Except for the superoptimal circulant matrix, IT of CNAS-GMRES with all the other circulant matrices are less than $10$, and the optimal empirical value of $\omega$ are almost the same, i.e., [0.21,1.01] ([0.41,0.81] for T. Chan's circulant matrix). These experiments show that most of the circulant matrices in the literature are efficient when they are applied to CNAS-GMRES, which makes us easy to choose the circulant approximation in CNAS-GMRES. The superoptimal circulant matrix performs poorly in our experiments, and a possible reason is that the optimal value $\omega$ of CNAS-GMRES in this case stays outside the search interval $(0,4]$.

\begin{table}[htbp]
	\setlength{\abovecaptionskip}{0pt}
	\setlength{\belowcaptionskip}{10pt} \centering{
		\caption{\label{tab:TNS_different_circulant_alp=1.9}
			IT of CNAS-GMRES with different circulant matrices, and the empirical optimal value of the parameter $\omega$ of CNAS-GMRES when $\alpha=1.9$, $M=6400, N=200$.}
		\begin{tabular}{llrccccccc}\specialrule{0em}{2pt}{2pt}\hline\specialrule{0em}{2pt}{2pt}
			
			Circulant Matrix & & &  IT & &  & & $\omega$  \\\specialrule{0em}{1pt}{1pt}\hline\specialrule{0em}{3pt}{3pt}
			
			T. Chan & &  & 7 & & & & 	[0.41,0.81]
			\\\specialrule{0em}{3pt}{3pt}
			
			Strang & &	 & 8 & &	 & &	[0.21,1.01]
			\\\specialrule{0em}{3pt}{3pt}

			R. Chan & &  &	8 & &	 & &	[0.21,1.01]
			\\\specialrule{0em}{3pt}{3pt}
			
			Modified Dirichlet kernel & &  & 8& &  & &	[0.21,1.01]
			\\\specialrule{0em}{3pt}{3pt}		
			
			von Hann kernel & &  &	8 & &	 & & [0.21,1.01]
			
			\\\specialrule{0em}{3pt}{3pt}
			
			Hamming kernel & &  & 8 & &	& &	[0.21,1.01]
			\\\specialrule{0em}{3pt}{3pt}
			
			Superoptimal & &  &	220 & &	 & &	[0.52,0.55]
			\\\specialrule{0em}{3pt}{3pt}
			\hline
	\end{tabular}}
\end{table}

\begin{figure}[htbp]
	\centering
	\subfloat{\includegraphics[scale=0.40]{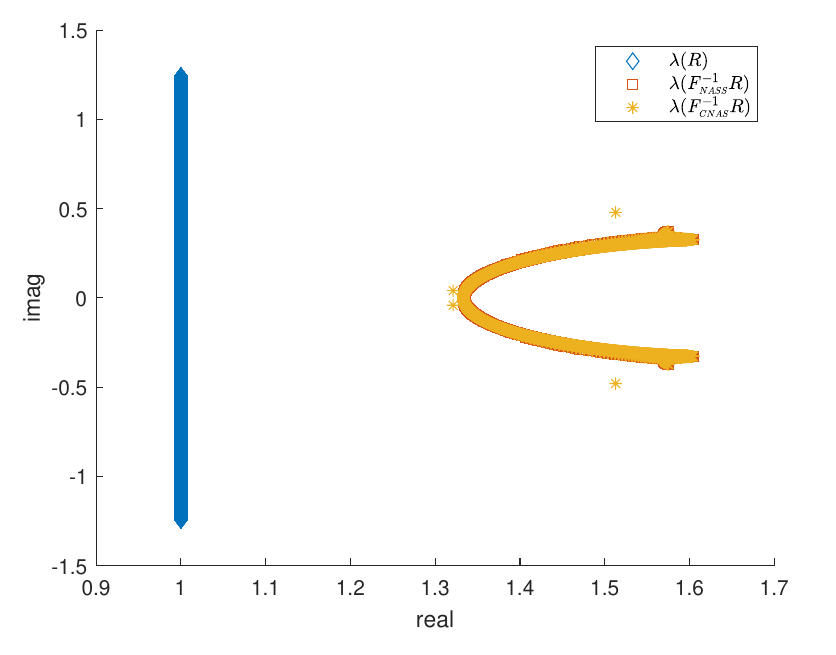}}
	\subfloat{\includegraphics[scale=0.40]{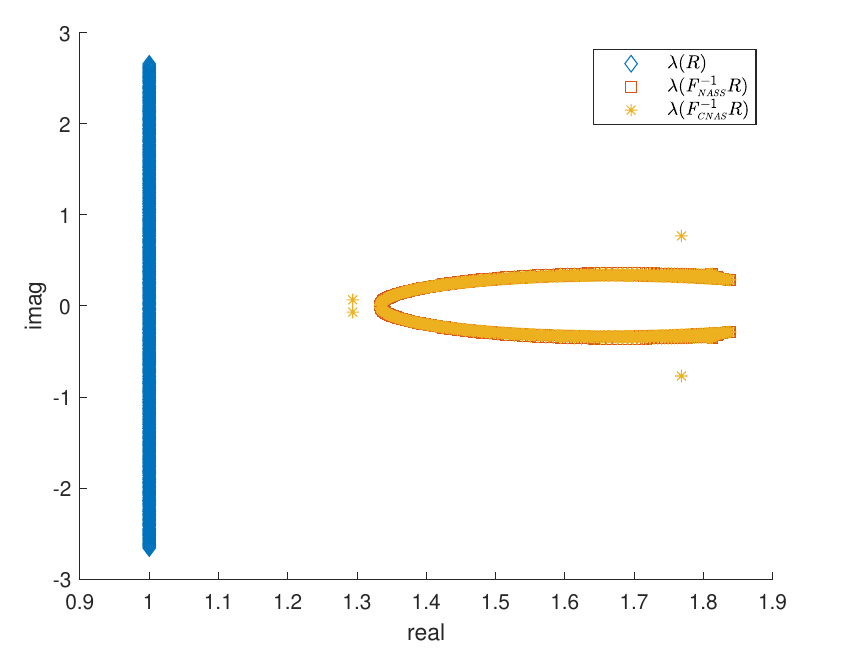}}
	
	\caption{
		The eigenvalue distribution of $\mathcal{R}$, $\mathcal{F}_{\mbox{\tiny NASS}}^{-1}\mathcal{R}$, and $\mathcal{F}_{\mbox{\tiny CNAS}}^{-1}\mathcal{R}$ when $\alpha=1.1$, $M=1600, N=200$ $(\text{left})$ and $M=3200, N=200$ $(\text{right})$.
	}
	\label{fig:eig1.1}
\end{figure}

\begin{figure}[htbp]
	\centering
	\subfloat{\includegraphics[scale=0.40]{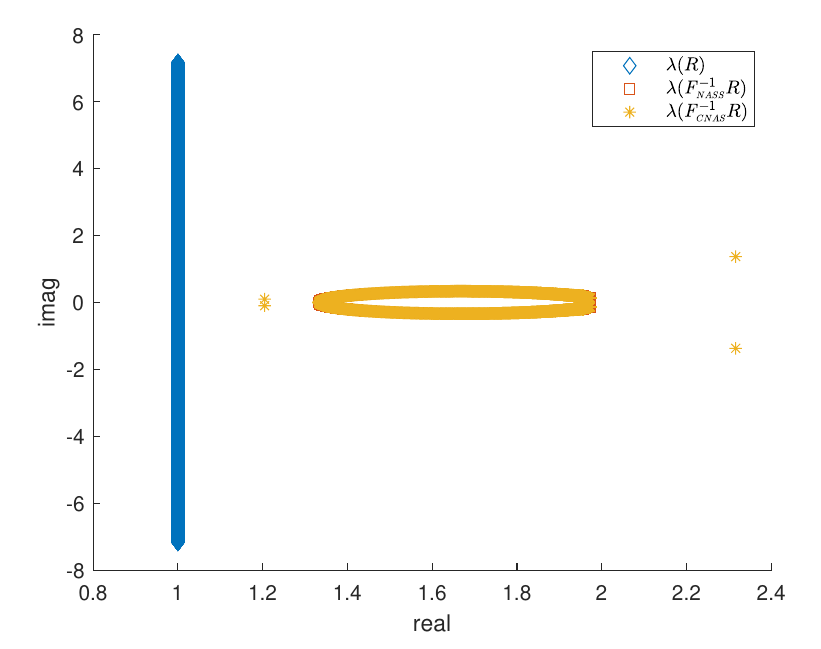}}
	\subfloat{\includegraphics[scale=0.40]{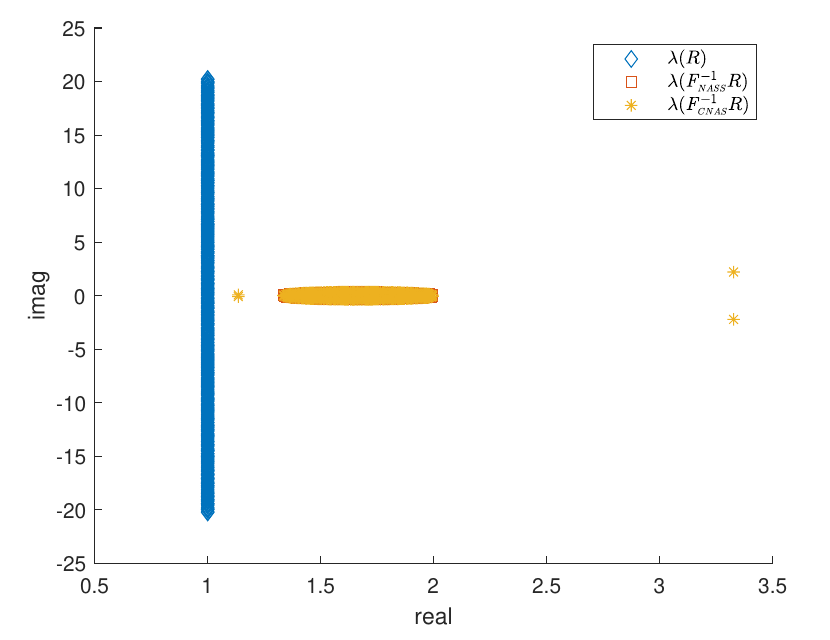}}

	\caption{
		The eigenvalue distribution of $\mathcal{R}$, $\mathcal{F}_{\mbox{\tiny NASS}}^{-1}\mathcal{R}$, and $\mathcal{F}_{\mbox{\tiny CNAS}}^{-1}\mathcal{R}$ when $\alpha=1.5$, $M=1600, N=200$ $(\text{left})$ and $M=3200, N=200$ $(\text{right})$.
	}
    \label{fig:eig1.5}
\end{figure}

\begin{figure}[htbp]
	\centering
	\subfloat{\includegraphics[scale=0.40]{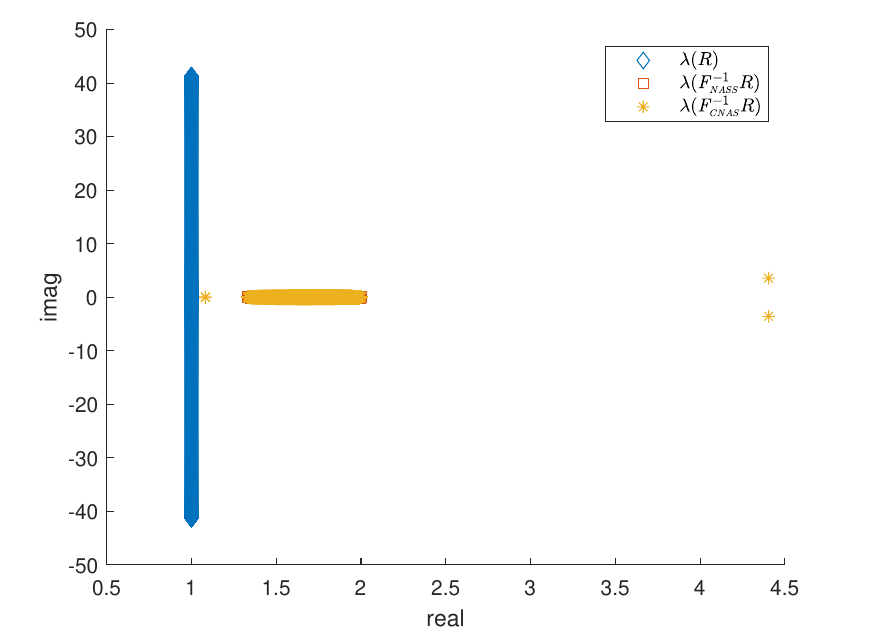}}
	\subfloat{\includegraphics[scale=0.40]{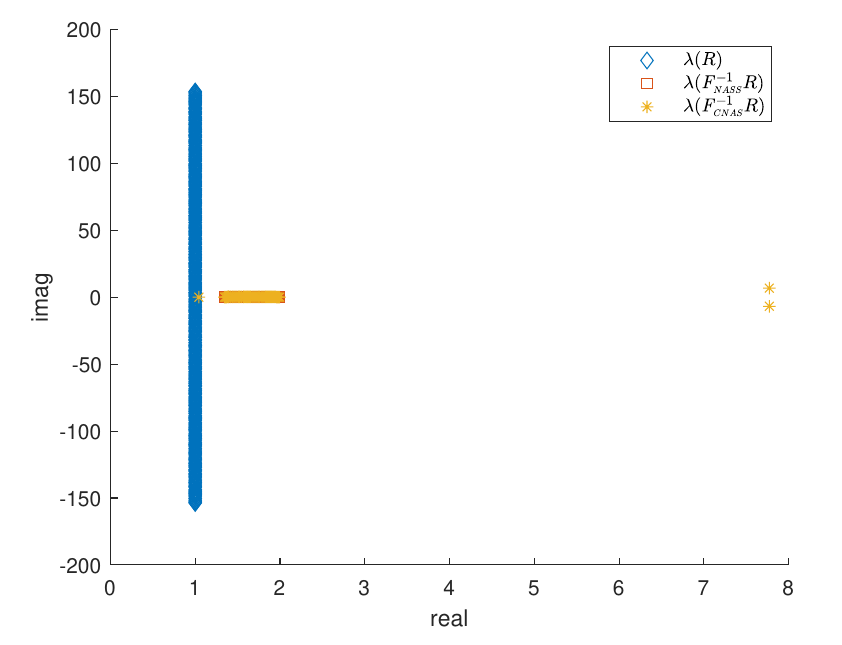}}

	\caption{
		The eigenvalue distribution of $\mathcal{R}$, $\mathcal{F}_{\mbox{\tiny NASS}}^{-1}\mathcal{R}$, and $\mathcal{F}_{\mbox{\tiny CNAS}}^{-1}\mathcal{R}$ when $\alpha=1.9$,  $M=1600, N=200$ $(\text{left})$ and $M=3200, N=200$ $(\text{right})$.
	}
	\label{fig:eig1.9}
\end{figure}

Figures \ref{fig:eig1.1}-\ref{fig:eig1.9} depict the eigenvalue distribution of the system matrix $\mathcal{R}$, the NASS preconditioned system matrix $\mathcal{F}_{\mbox{\tiny NASS}}^{-1}\mathcal{R}$, and the CNAS preconditioned system matrix $\mathcal{F}_{\mbox{\tiny CNAS}}^{-1}\mathcal{R}$ when $\alpha=1.1:0.4:1.9$, $M=1600$, $3200$. In Figures \ref{fig:eig1.1}-\ref{fig:eig1.9}, the left plots are related to the case $M=1600$, and the right plots are related to the case $M=3200$.
In Figure \ref{fig:eig1.1}, the real parts of the eigenvalues of $\mathcal{R}$ are 1, and the imaginary parts are distributed from $-1.3$ to $1.3$ when $M=1600$, and from $-2.7$ to $2.7$ when $M=3200$.
In Figure \ref{fig:eig1.5}, the real parts of the eigenvalues of $\mathcal{R}$ are 1, and the imaginary parts are distributed from $-8$ to $8$ when $M=1600$, and from $-20$ to $20$ when $M=3200$.
In Figure \ref{fig:eig1.9}, the real parts of the eigenvalues of $\mathcal{R}$ are 1, and the imaginary parts are distributed from $-42$ to $42$ when $M=1600$, and from $-150$ to $150$ when $M=3200$. Meanwhile, the real parts of the eigenvalues of $\mathcal{F}_{\mbox{\tiny NASS}}^{-1}\mathcal{R}$ and $\mathcal{F}_{\mbox{\tiny CNAS}}^{-1}\mathcal{R}$ are distributed from $1.3$ to $2$, and the imaginary parts are distributed from  $-0.5$ to $0.5$ in all the plots of Figures \ref{fig:eig1.1}-\ref{fig:eig1.9}.
Obviously, the eigenvalues of $\mathcal{F}_{\mbox{\tiny NASS}}^{-1}\mathcal{R}$ and $\mathcal{F}_{\mbox{\tiny CNAS}}^{-1}\mathcal{R}$ are more clustered than those of $\mathcal{R}$ (especially for larger $\alpha$). In addition, most of the eigenvalues of $\mathcal{F}_{\mbox{\tiny CNAS}}^{-1}\mathcal{R}$ are clustered around those of $\mathcal{F}_{\mbox{\tiny NASS}}^{-1}\mathcal{R}$, and only several eigenvalues of $\mathcal{F}_{\mbox{\tiny CNAS}}^{-1}\mathcal{R}$ drift a little further away from those of $\mathcal{F}_{\mbox{\tiny NASS}}^{-1}\mathcal{R}$, which fits the prediction of Remarks \ref{rmk:Preconditioning1}-\ref{rmk:Preconditioning2}. The eigenvalue distribution of the preconditioned cases keeps almost the same when $M$ increases from $1600$ to $3200$, which indicates the space mesh size independent convergence property of NASS-GMRES and CNAS-GMRES.

\begin{figure}[htbp]
	\centering
	\subfloat{\includegraphics[scale=0.4]{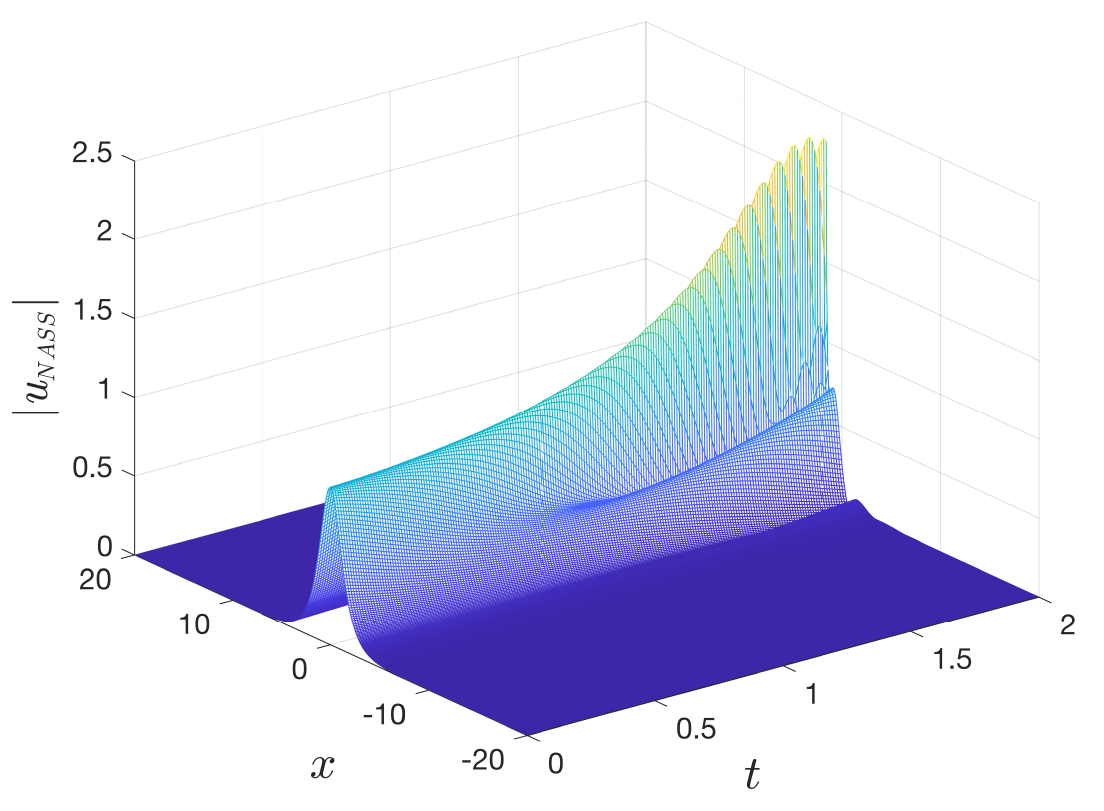}}
	\subfloat{\includegraphics[scale=0.4]{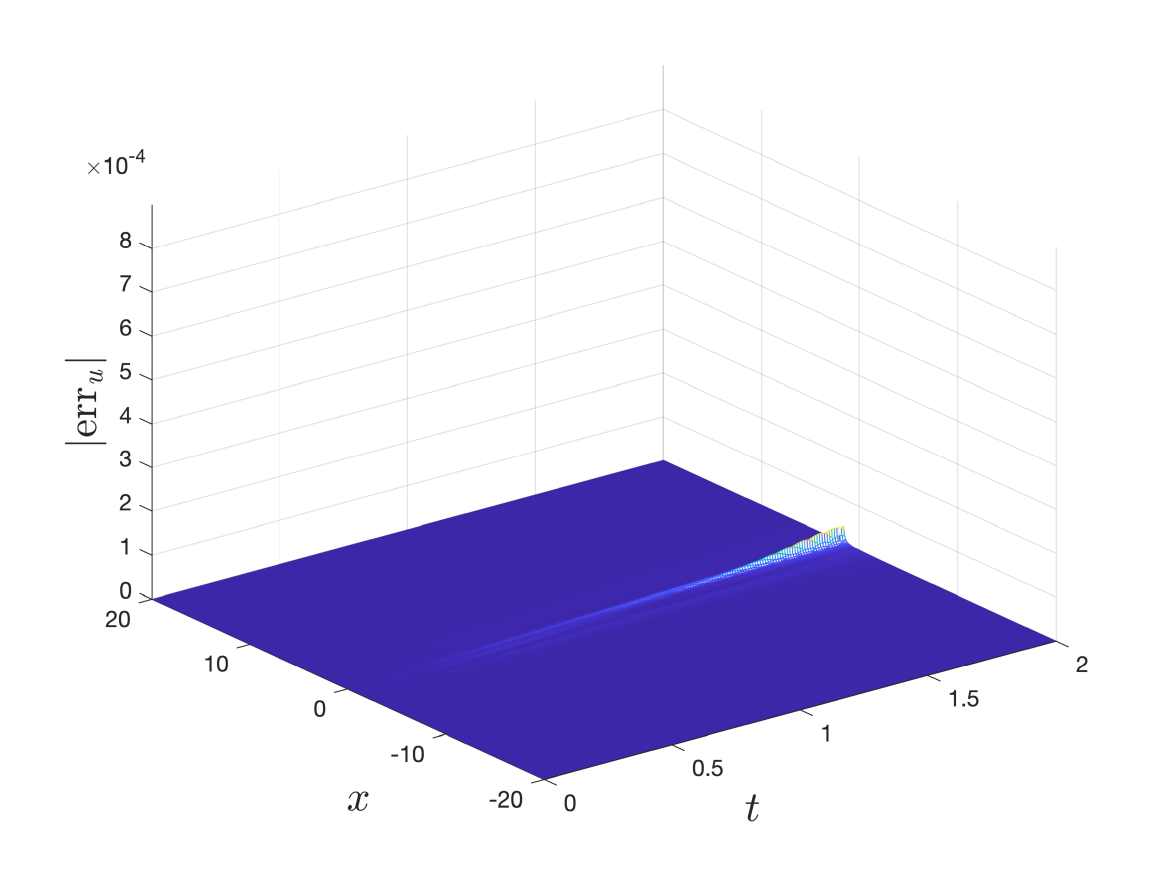}}
	\caption{
		The numerical solution (left) and its error (right) of  the space fractional DNLS equation (\ref{equ:DNLS}) when $\alpha=1.1$, $M=800, N=200$.}
	\label{fig:singal_1.1}
\end{figure}

\begin{figure}[htbp]
	\centering
	\subfloat{\includegraphics[scale=0.40]{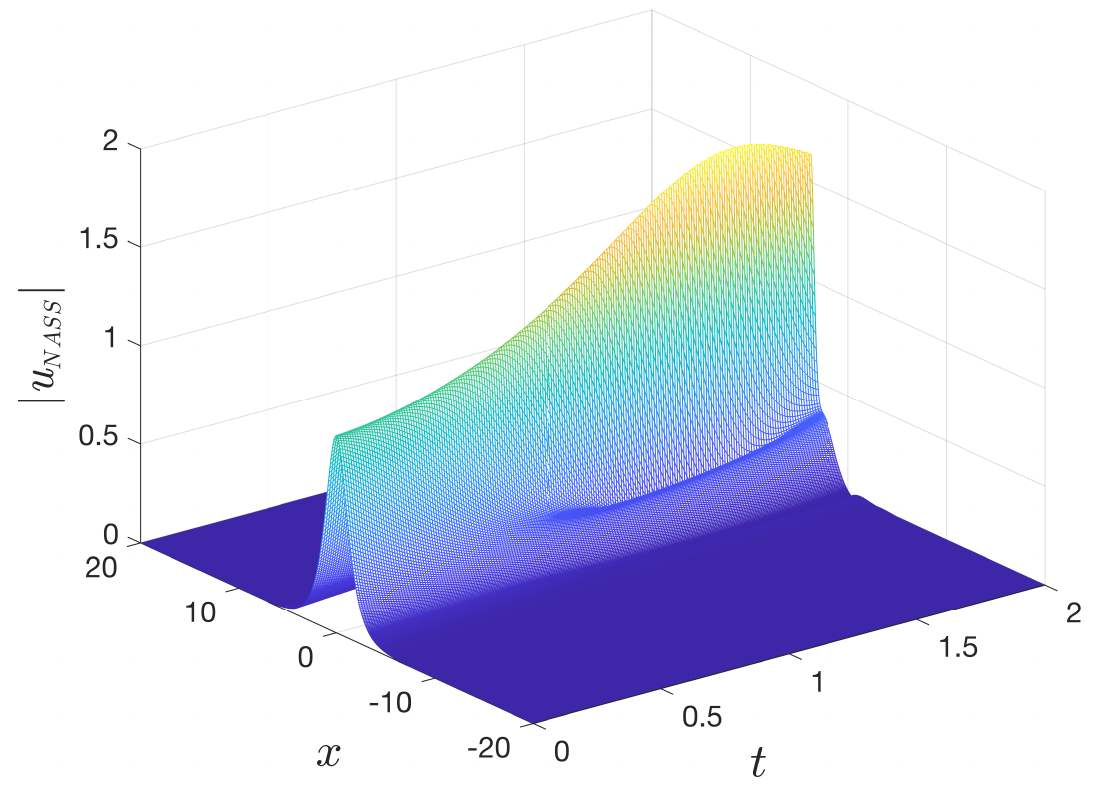}}
	\subfloat{\includegraphics[scale=0.40]{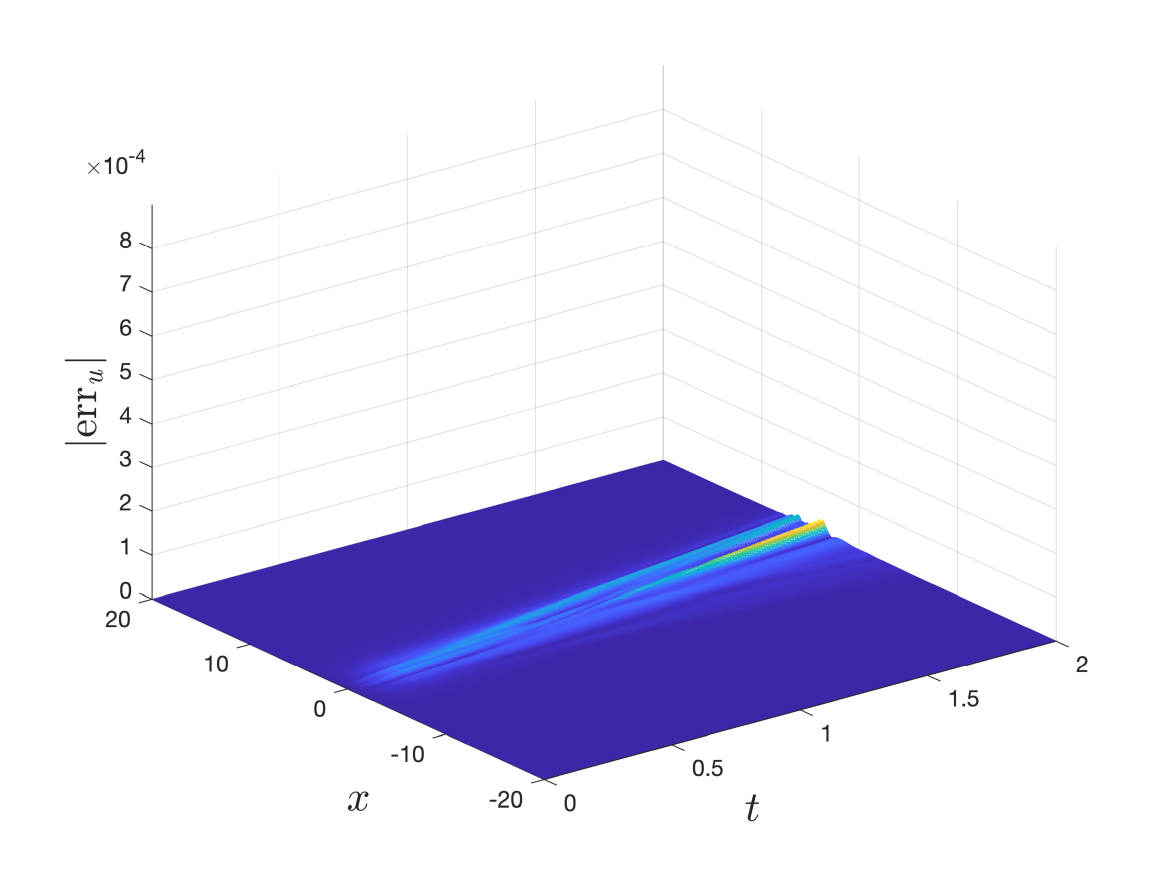}}
	\caption{
		The numerical solution (left) and its error (right) of the space fractional DNLS equation (\ref{equ:DNLS}) when $\alpha=1.5$, $M=800, N=200$.}
	\label{fig:singal_1.5}
\end{figure}

\begin{figure}[htbp]
	\centering
	\subfloat{\includegraphics[scale=0.4]{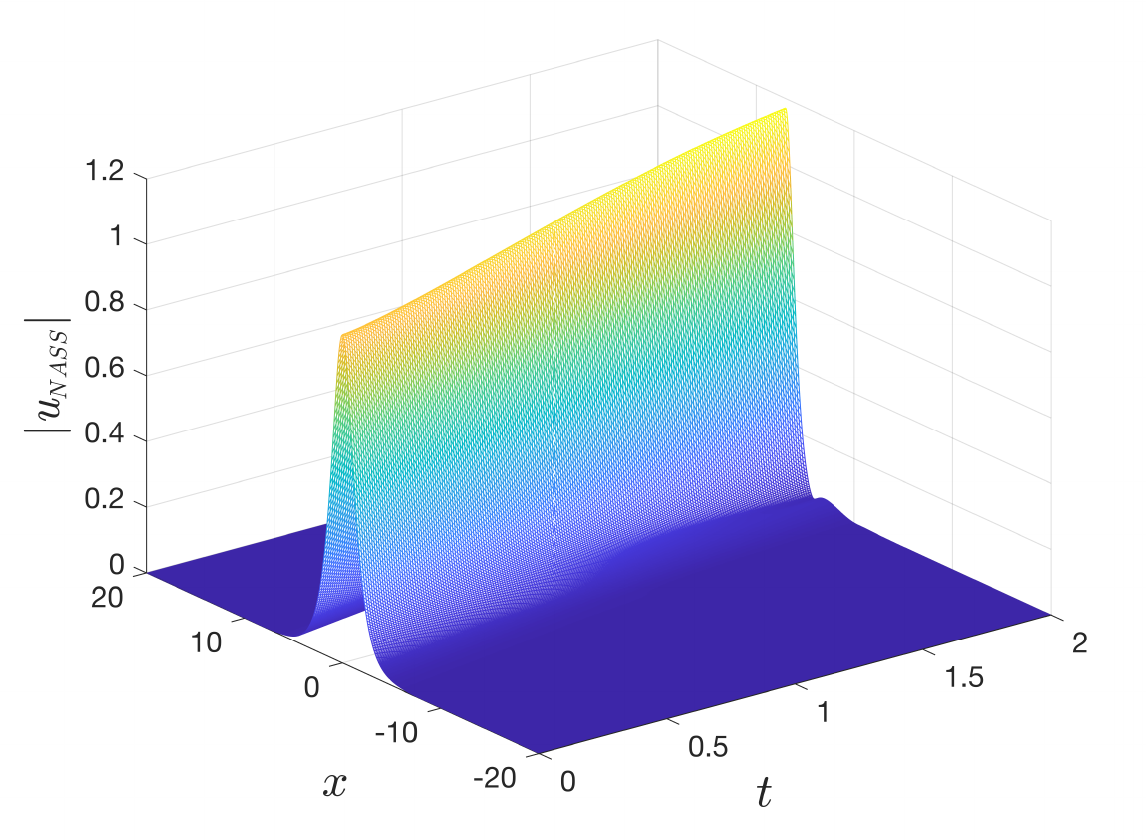}}
	\subfloat{\includegraphics[scale=0.4]{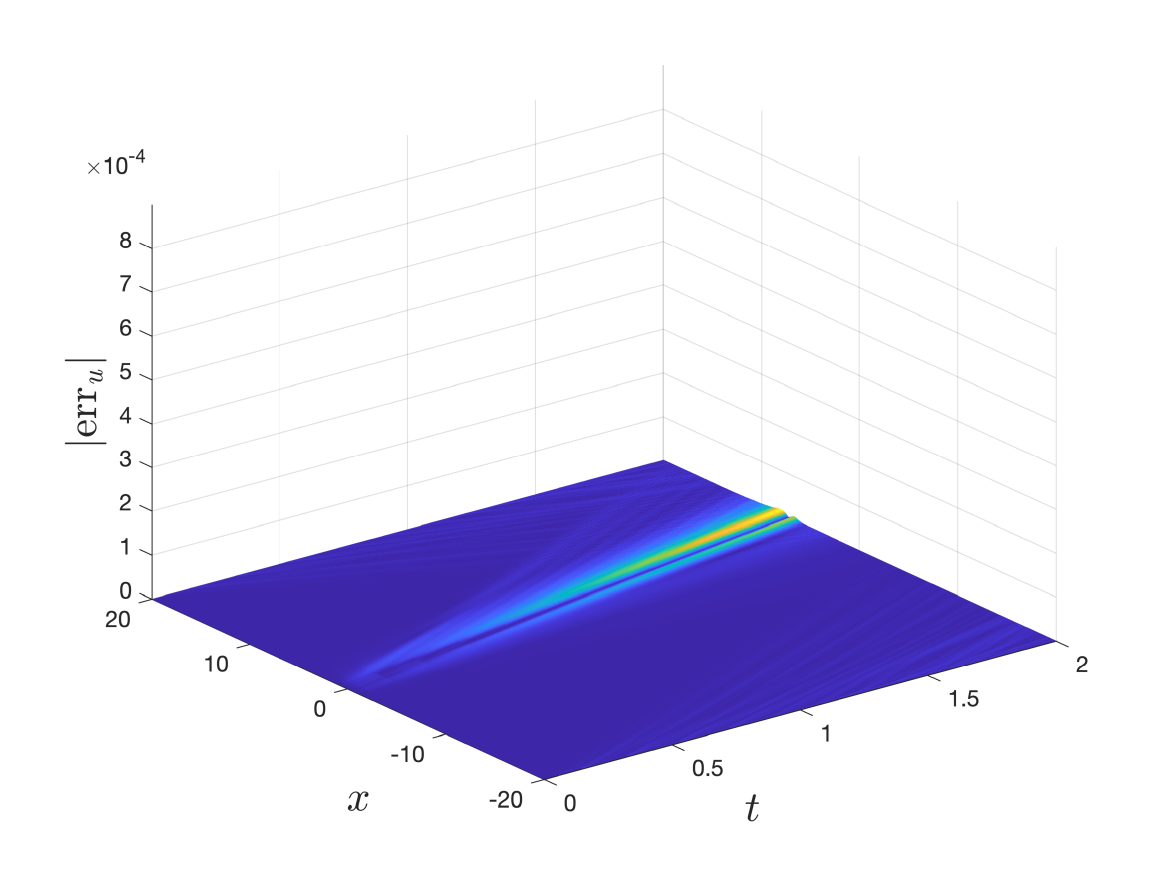}}
	\caption{
		The numerical solution (left) and its error (right)  of  the space fractional DNLS equation (\ref{equ:DNLS}) when $\alpha=1.9$, $M=800, N=200$.}
	\label{fig:singal_1.9}
\end{figure}

\begin{figure}[htbp]
	\centering
	\subfloat{\includegraphics[scale=0.4]{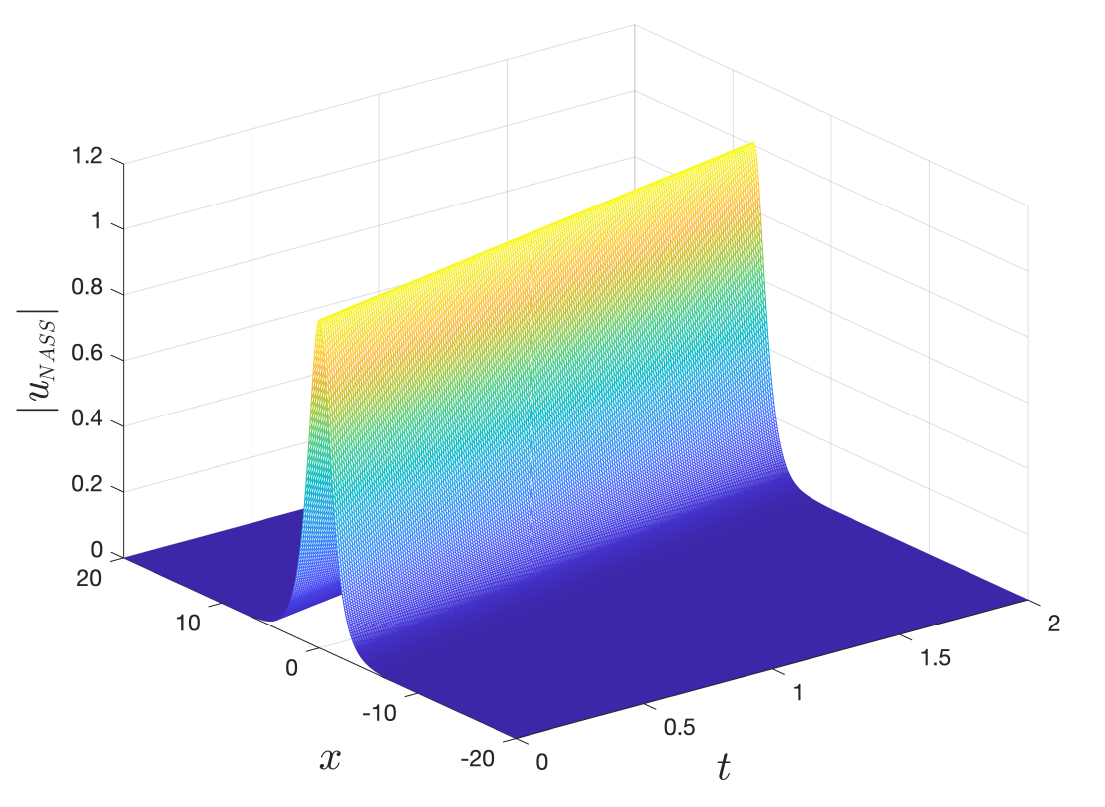}}
	\subfloat{\includegraphics[scale=0.4]{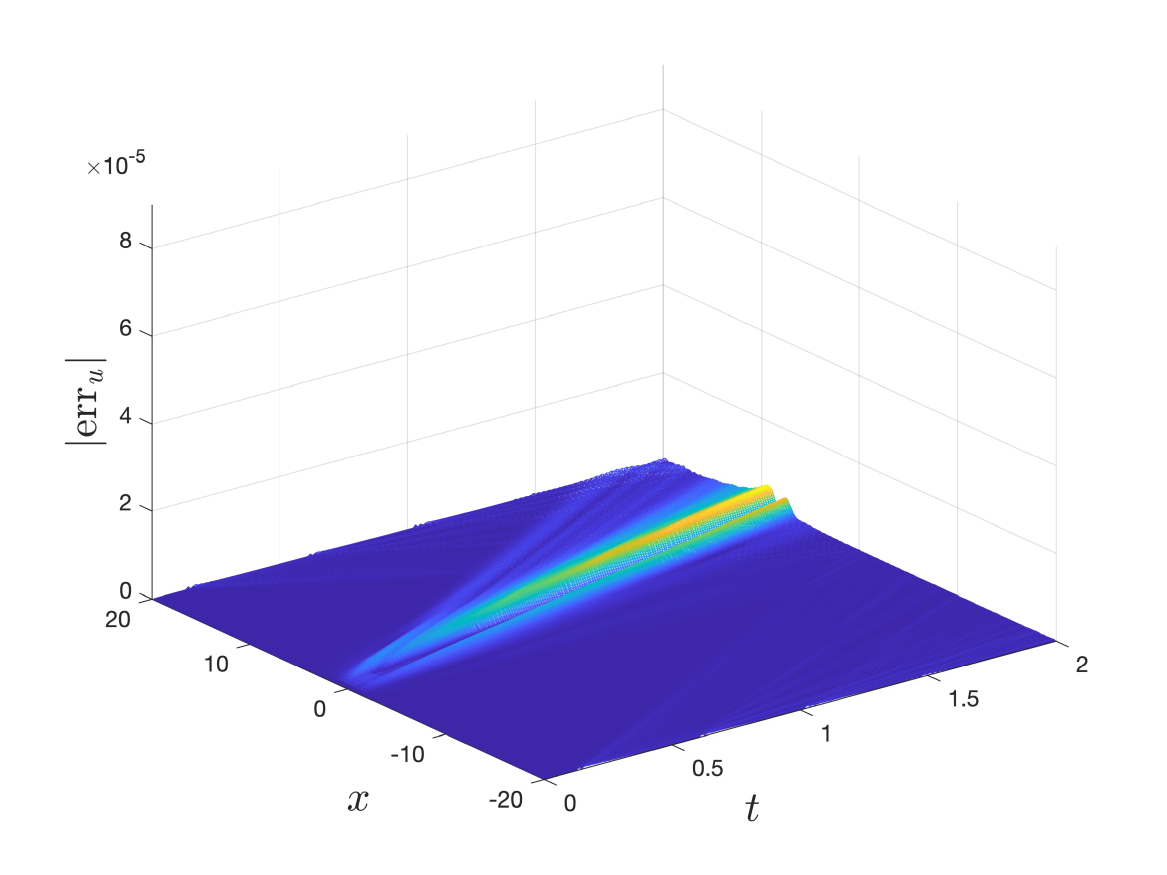}}
	\caption{
		The numerical solution (left) and its error (right) of the space fractional DNLS equation (\ref{equ:DNLS})  when $\alpha=2$, $M=800, N=200$.}
	\label{fig:singal_2}
\end{figure}

Figures \ref{fig:singal_1.1}-\ref{fig:singal_2} depict the numerical solution (left), i.e., $u_{\text{\tiny{CNAS}}}$, obtained by  CNAS-GMRES, and its error (right), i.e., $\text{err}_{u} = |u_{\text{\tiny CNAS}}-u_{\text{\tiny GE}}|$ with the exact solution $u_{\text{\tiny GE}}$ of the LICD scheme (i.e., the solution $u_{\text{\tiny GE}}$ obtained by GE) for the space fractional DNLS equations (\ref{equ:DNLS}) when $M=800, N=200$, $\alpha=1.1:0.4:1.9$ and $\alpha=2$.
It is  observed that the fractional order $\alpha$ will affect the shape of the wave front both in height and width. When $\alpha$ gets smaller, the shape of the wave front will change more quickly. When $\alpha$ tends to 2, the wave front  of the space fractional DNLS equation converges to the wave front of the non-fractional one, and the height of the wave front is stable.		
 In addition, the error between the numerical solution and the exact solution of the LICD scheme remains as small as around $10^{-4}$ in the whole computational space-time domain, which shows that the numerical solution obtained by CNAS-GMRES is reliable.

\begin{figure}[htbp]
	\centering
	\includegraphics[scale=0.6]{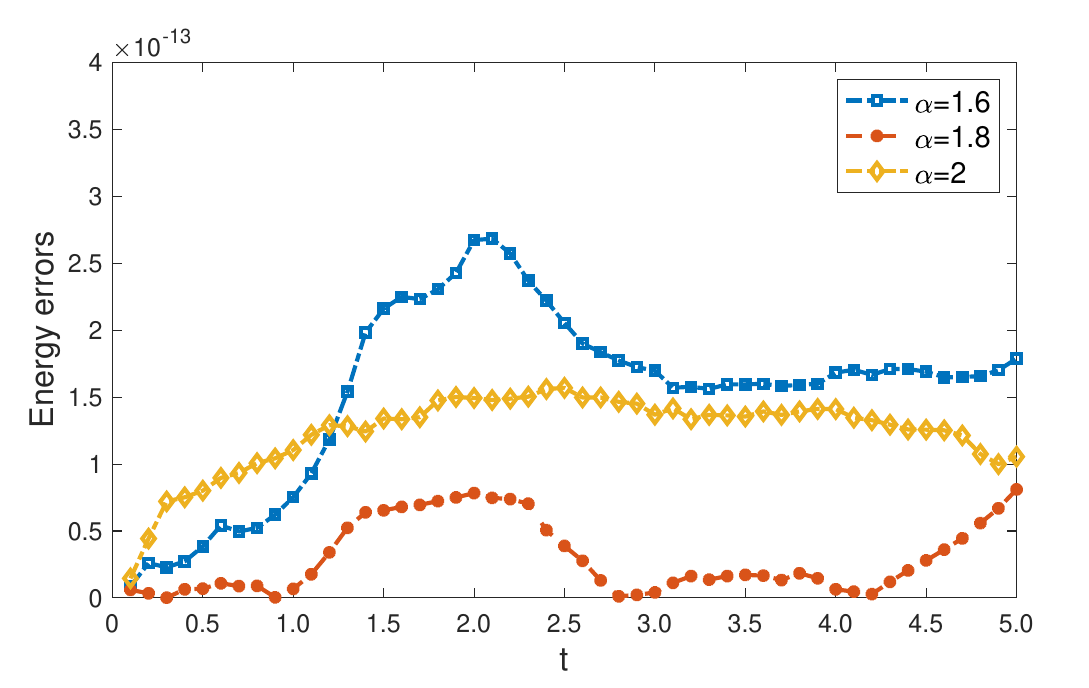}
	\caption{
		The relative errors of the discrete energy, i.e., $|(E^n-E^0)/E^0|$, when $h=0.2$, $\tau = 0.001$.
	}
	\label{fig:energy1}
\end{figure}

Table \ref{mass1} and Figure \ref{fig:energy1} report the relative errors of the discrete mass and energy  for different values of $\alpha$ respectively. The space mesh size and time step size are $h=0.2$ and $\tau=0.05$ for Table \ref{mass1}, and $h=0.2$ and $\tau=0.001$ for Figure \ref{fig:energy1}. At each time level, the block linear system is solved by CNAS-GMRES, and terminated when the $\ell_2$-norm relative residual of the block linear system reduced below $10^{-15}$. The errors reported in Table \ref{mass1} and Figure \ref{fig:energy1} stay very small, which means that the linear solver CNAS-GMRES maintains the conservation properties of the LICD scheme.

\begin{table}[htbp]
	\renewcommand\arraystretch{1.6}
	
	\tabcolsep=0.6cm
	\setlength{\abovecaptionskip}{2pt}
	\setlength{\belowcaptionskip}{16pt} \centering{
		\caption{
			\label{mass1}
			The relative errors of the discrete mass, i.e., $|(Q^{n}-Q^{0})/Q^{0}|$, when $h=0.2,\tau=0.05.$}
\begin{tabular}{lllll}
	
	\hline & $t=1$ & $t=2$ & $t=3$ & $t=4$ \\
	\hline$\alpha=1.4$ & $4.8850 \mathrm{e}-015$ & $7.5495 \mathrm{e}-015$ & $6.4393 \mathrm{e}-015$ & $9.1038 \mathrm{e}-015$ \\
	
	$\alpha=1.7$ & $5.9952\mathrm{e}-015$ & $4.2188 \mathrm{e}-015$ & $3.3307\mathrm{e}-015$ & $3.7748\mathrm{e}-015$ \\
	
	$\alpha=1.9$ & $1.1102 \mathrm{e}-015$ & $2.2209\mathrm{e}-015$ & $3.1086 \mathrm{e}-015$ &  $1.5543 \mathrm{e}-015$\\
	
	$\alpha=2$ & $6.6613 \mathrm{e}-016$ & $5.7732 \mathrm{e}-015$ & $5.3291 \mathrm{e}-015$ & $6.6615 \mathrm{e}-015$\\
	\hline
\end{tabular}}
\end{table}

\subsection{The CNLS case}

The experiments here are carried out on the coupled case, i.e., the truncated system of CNLS equations
\begin{align}\label{equ:CNLS}
	\begin{cases}
		\imath u_t-\gamma(-\Delta)^{\frac{\alpha}{2}}u+\rho(\vert u \vert^2+\beta\vert v \vert^2)u=0, \qquad  \\
		\imath v_t-\gamma(-\Delta)^{\frac{\alpha}{2}}v+\rho(\vert v \vert^2+\beta\vert u \vert^2)v=0, \qquad
	\end{cases}
	-20\le x \le 20,	\quad 0 < t \le \mbox{T},
\end{align}
subjected to the initial and boundary conditions	
\begin{align}\label{equ:CNLS_InitBdry}
	\begin{cases}
		u(x,0)=\text{sech}(x+5) \ e^{3\imath x}, \quad v(x,0)=\text{sech}(x-5) \ e^{-3\imath x},\\
		u(-20,t)=u(20,t)=0, \quad v(-20,t)=v(20,t)=0,
	\end{cases}
\end{align}
where $\gamma=1$, $\rho=1$, $\beta=1$, $1<\alpha\le 2$. The LICD scheme applied to (\ref{equ:CNLS}) and (\ref{equ:CNLS_InitBdry}) leads to the discretized space fractional CNLS equations. It requires to solve two complex symmetric linear systems of the form (\ref{equ3}) successively at each time level $t_n$ for $1 < n \le N$, which is equivalent to solve two block linear systems of the form (\ref{positiveBlockForm}). In the sequel, `CPU' and `IT' are the total computing time in seconds and the total iteration counts for solving two coupled linear systems at $t_{n}$.

\begin{figure}[htbp]
	\centering
	\includegraphics[scale=0.38]{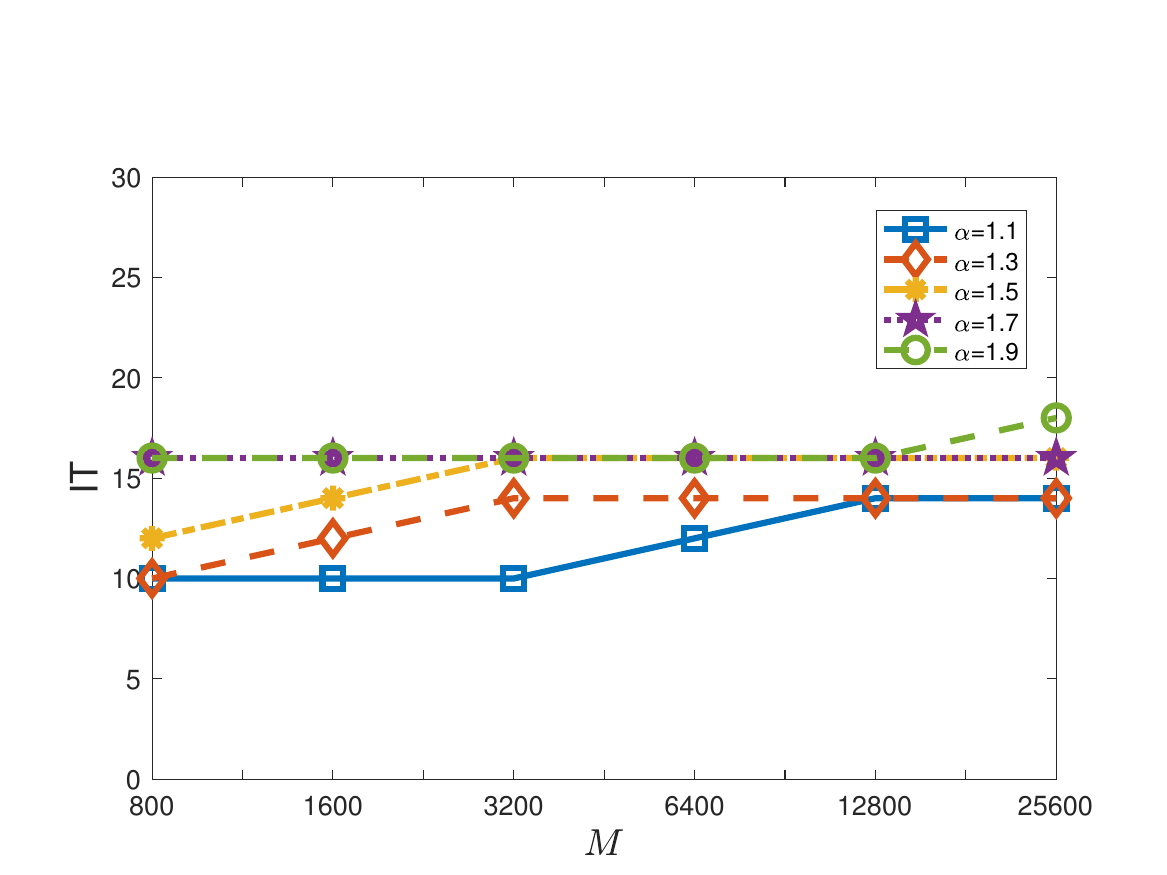}
	\caption{
		The curves of IT of CNAS-GMRES versus the space mesh size $M$ when $\alpha=1.1:0.2:1.9$, $N=200$.
	}
	\label{CMit}
\end{figure}

Figure \ref{CMit} depicts the curves of IT of CNAS-GMRES versus the space mesh size $M$ of the LICD scheme applied to the space fractional CNLS equations when $\alpha=1.1:0.2:1.9$, $M=800$, $1600$, $3200$, $6400$, $12800$, $25600$ and $N=200$. The empirical optimal value of $\omega$ is selected. The result is similar to the case of the discretized space fractional DNLS equations, which confirms that  CNAS-GMRES is also space mesh size independent in the coupled case.

\begin{figure}[htbp]
	\centering
	\subfloat{\includegraphics[scale=0.276]{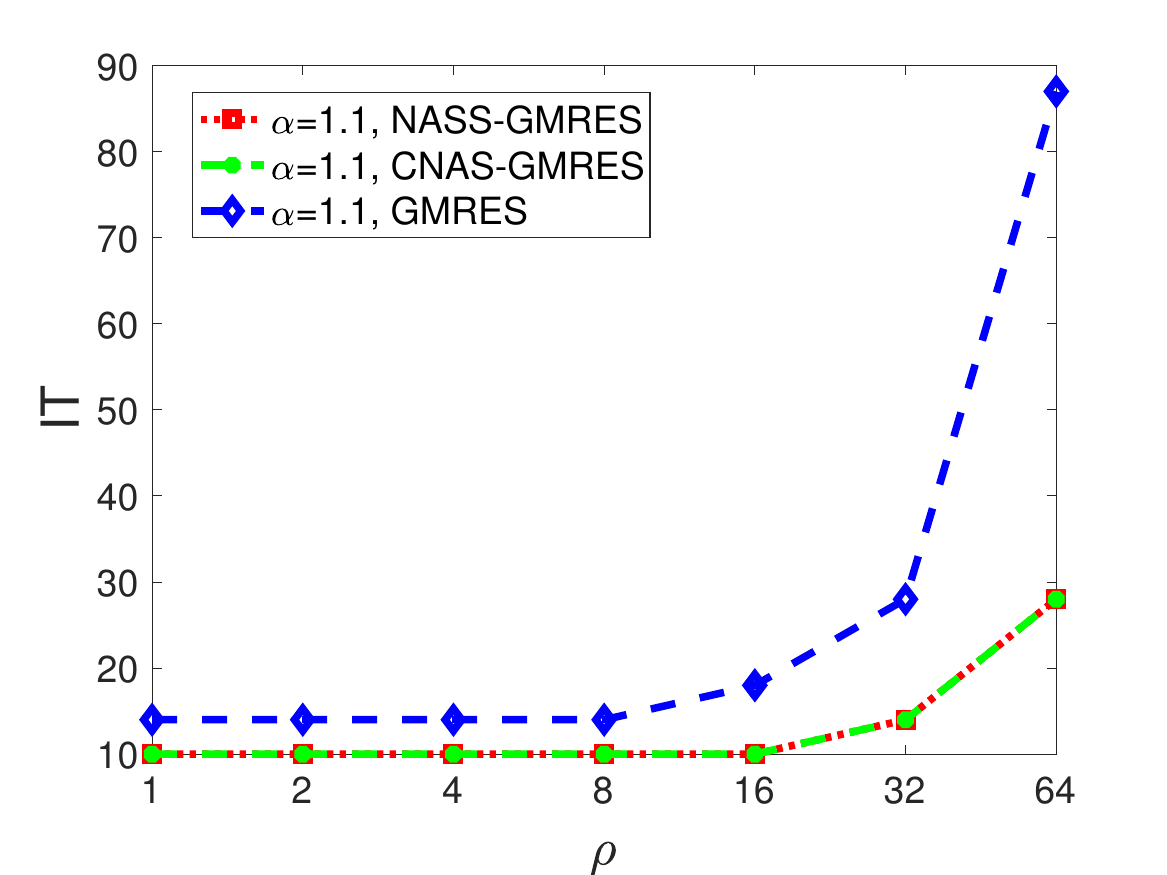}}
	\subfloat{\includegraphics[scale=0.276]{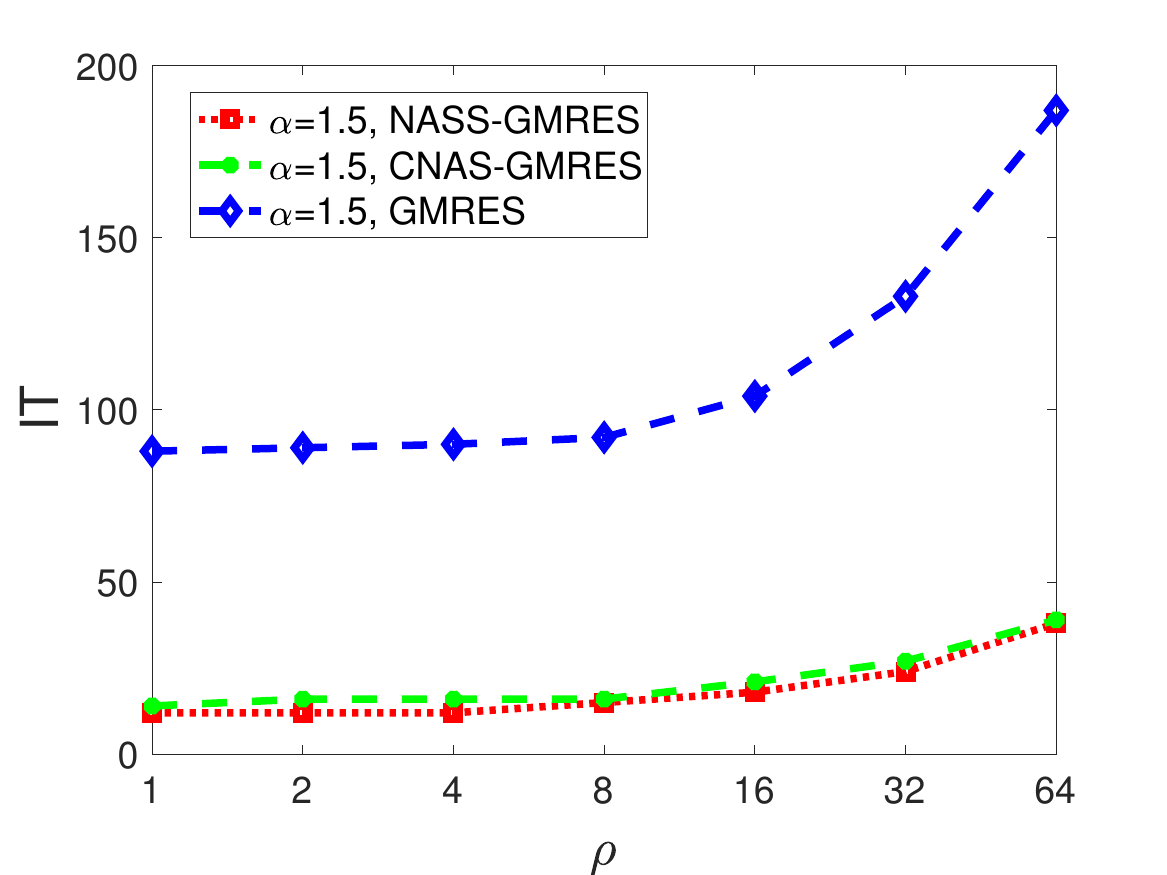}}
	\subfloat{\includegraphics[scale=0.276]{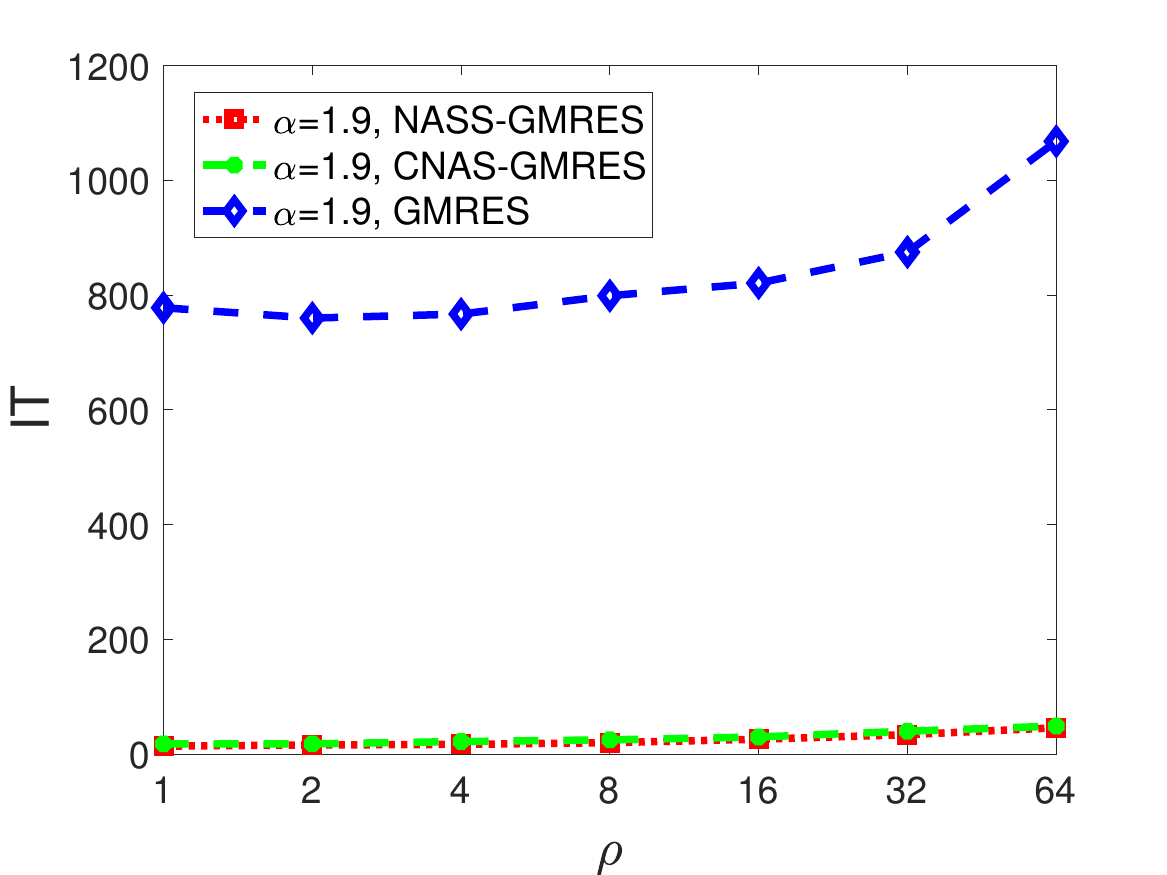}}
	\caption{
	The curves of IT of NASS-GMRES, CNAS-GMRES and GMRES versus the nonlinear term parameter $\rho$ when $\alpha=1.1:0.4:1.9$, $M=1600, N=200$.
	}
	\label{rhoit}
\end{figure}

Figure \ref{rhoit} describes the influence of the strength of the nonlinear term (controlled by the parameter $\rho$) on the iteration counts of GMRES, NASS-GMRES and CNAS-GMRES when $\alpha=1.1:0.4:1.9 $, $M=1600$. The parameter $\rho$ grows from 1 to 64. In these plots, IT of NASS-GMRES and CNAS-GMRES is very small, and increases slowly as $\rho$ grows. Meanwhile, IT of GMRES is very large, and increases quickly. It can be seen that the stronger the nonlinear term is, the more difficult it is to solve the coupled linear systems, and the new preconditioned GMRES methods can effectively deal with the strongly nonlinear case.

Tables \ref{tab:alp=1.1_beta=1}-\ref{tab:alp=1.9_beta=1} provide CPU and IT of GMRES, CNAS-GMRES, and GE when $\alpha=1.1:0.2:1.9$, $M=3200$, $6400$, $12800$, $25600$, $N=200$. In these tables, `--' represents that GMRES is not convergent in a prescribed number of iterations, `N/A' represents that the data of IT for GE is not available. The empirical optimal parameters of CNAS-GMRES corresponding to the results in Tables \ref{tab:alp=1.1_beta=1}-\ref{tab:alp=1.9_beta=1} are listed in Table \ref{tab:omega alp=1.1}. Specifically, `$\omega_u$' and `$\omega_v$' are denoted by the empirical optimal parameters of CNAS-GMRES for the block linear systems (\ref{positiveBlockForm}) related to $u$ and $v$, respectively. According to Tables \ref{tab:alp=1.1_beta=1}-\ref{tab:alp=1.9_beta=1}, GE requires the largest CPU in all the tests, and CPU of  CNAS-GMRES is always smaller than CPU of GMRES. Besides, IT of GMRES  increases quickly when the fractional order $\alpha$ and the space mesh size $M$ increase, indicating that the system (\ref{equ:CNLS}) is more difficult to solve as $\alpha$ grows. Meanwhile, IT of CNAS-GMRES is smaller than IT of GMRES, especially for larger $\alpha$ and $M$. In summary, CNAS-GMRES improves the computational eﬀiciency  significantly, and behaves independent on the space mesh size $M$ and almost insensitive to the fractional order $\alpha$.

\begin{table}[htbp]
	\setlength{\abovecaptionskip}{0pt}
	\setlength{\belowcaptionskip}{10pt} \centering{
		\caption{\label{tab:alp=1.1_beta=1} 
			CPU and IT of GMRES,  CNAS-GMRES, and GE when $\alpha=1.1$, $N=200$.}
		\begin{tabular}{lcccccccc}\specialrule{0em}{2pt}{2pt}\hline\specialrule{0em}{2pt}{2pt}
			$M$  & 3200 &  & 6400  &  & 12800  &  &  25600  & \\\cmidrule(l){2-3}\cmidrule(l){4-5}\cmidrule(l){6-7}\cmidrule(l){8-9}
			&  CPU &  IT & CPU &  IT  & CPU &  IT& CPU &  IT   \\\specialrule{0em}{1pt}{1pt}\hline\specialrule{0em}{3pt}{3pt}
			
			GMRES & 	8.33E-02 &	22 &	5.69E-01 &	50 & 2.67E-00 & 110 & 1.10E+01	& 236
			\\\specialrule{0em}{3pt}{3pt}

			CNAS-GMRES & 	6.10E-02 &	10&	2.61E-01 &	12 & 6.98E-01 & 14 & 1.21E-00 &	14			
			\\\specialrule{0em}{3pt}{3pt}
			
			GE  & 2.10E+01	 &	N/A &	1.67E+02 &	N/A & 1.36E+03  & N/A & 1.88E+04	&  N/A
			\\\specialrule{0em}{3pt}{3pt}	\hline
			
	\end{tabular}}
\end{table}

\begin{table}[htbp]
	\setlength{\abovecaptionskip}{0pt}
	\setlength{\belowcaptionskip}{10pt} \centering{
		\caption{\label{tab:alp=1.3_beta=1} 
			CPU and IT of GMRES, CNAS-GMRES, and GE when $\alpha=1.3$, $N=200$.}
		\begin{tabular}{lcccccccc}\specialrule{0em}{2pt}{2pt}\hline\specialrule{0em}{2pt}{2pt}
			$M$  & 3200 &  & 6400  &  & 12800  &  &  25600  & \\\cmidrule(l){2-3}\cmidrule(l){4-5}\cmidrule(l){6-7}\cmidrule(l){8-9}
			&  CPU &  IT & CPU &  IT  & CPU &  IT& CPU &  IT   \\\specialrule{0em}{1pt}{1pt}\hline\specialrule{0em}{3pt}{3pt}
			
			GMRES & 	3.84E-01 &	66 &	2.82E-00 &	182 & 2.47E+01 & 514 & 1.61E+02	& 1087
			\\\specialrule{0em}{3pt}{3pt}

			CNAS-GMRES & 	6.52E-02 &	14&	2.53E-01 &	14 & 7.12E-01 & 14 & 1.23E-00 &	14			
			\\\specialrule{0em}{3pt}{3pt}
			
			GE & 2.10E+01 &N/A	 &	1.64E+02 &	N/A&	1.36E+03 & N/A &  1.94E+04& 	N/A
			\\\specialrule{0em}{3pt}{3pt}	\hline
			
	\end{tabular}}
\end{table}

\begin{table}[htbp]
	\setlength{\abovecaptionskip}{0pt}
	\setlength{\belowcaptionskip}{10pt} \centering{
		\caption{\label{tab:alp=1.5_beta=1} 
			CPU and IT of GMRES, CNAS-GMRES, and GE when $\alpha=1.5$, $N=200$.}
		\begin{tabular}{lcccccccc}\specialrule{0em}{2pt}{2pt}\hline\specialrule{0em}{2pt}{2pt}
			$M$  & 3200 &  & 6400  &  & 12800  &  &  25600  &   \\\cmidrule(l){2-3}\cmidrule(l){4-5}\cmidrule(l){6-7}\cmidrule(l){8-9}
			&  CPU &  IT & CPU &  IT  & CPU &  IT& CPU &  IT   \\\specialrule{0em}{1pt}{1pt}\hline\specialrule{0em}{3pt}{3pt}
			
			GMRES & 	1.86E-00 &	256 &	2.89E+01 &	832 & 3.26E+02 & 2700 & --	& --
			\\\specialrule{0em}{3pt}{3pt}
			
			CNAS-GMRES & 	7.24E-02 &	16&	2.60E-01 &	16& 7.56E-01 & 16 & 1.26E-00 &	16			
			\\\specialrule{0em}{3pt}{3pt}
			
		   GE & 2.10E+01	 &N/A &	1.66E+02 &N/A	 & 1.34e+03 &N/A & 1.92E+04	& N/A
			\\\specialrule{0em}{3pt}{3pt}	\hline
			
	\end{tabular}}
\end{table}

\begin{table}[htbp]
	\setlength{\abovecaptionskip}{0pt}
	\setlength{\belowcaptionskip}{10pt} \centering{
		\caption{\label{tab:alp=1.7_beta=1} 
			CPU and IT of GMRES, CNAS-GMRES, and GE when $\alpha=1.7$, $N=200$.}
		\begin{tabular}{lcccccccc}\specialrule{0em}{2pt}{2pt}\hline\specialrule{0em}{2pt}{2pt}
			$M$  & 3200 &  & 6400  &  & 12800  &  &  25600  & \\\cmidrule(l){2-3}\cmidrule(l){4-5}\cmidrule(l){6-7}\cmidrule(l){8-9}
			&  CPU &  IT & CPU &  IT  & CPU &  IT& CPU &  IT   \\\specialrule{0em}{1pt}{1pt}\hline\specialrule{0em}{3pt}{3pt}
			
			GMRES & 	1.31E+01 &	1086 &	4.46E+02 & 4056  & -- & -- & --	& --
			\\\specialrule{0em}{3pt}{3pt}
			
			CNAS-GMRES & 	9.40E-02 &	16&	2.95E-01 &	16 & 8.09E-01 & 16 & 1.53E-00 &	16			
			\\\specialrule{0em}{3pt}{3pt}
			
			GE & 2.11E+01	 &	N/A &	1.68E+02 &	N/A & 1.39E+03  &N/A & 	1.96E+04 & N/A
			\\\specialrule{0em}{3pt}{3pt}	\hline
			
	\end{tabular}}
\end{table}

\begin{table}[htbp]
	\setlength{\abovecaptionskip}{0pt}
	\setlength{\belowcaptionskip}{10pt} \centering{
		\caption{\label{tab:alp=1.9_beta=1} 
			CPU and IT of GMRES, CNAS-GMRES, and GE when $\alpha=1.9$, $N=200$.}
		\begin{tabular}{lcccccccc}\specialrule{0em}{2pt}{2pt}\hline\specialrule{0em}{2pt}{2pt}
			$M$  & 3200 &  & 6400  &  & 12800  &  &  25600  & \\\cmidrule(l){2-3}\cmidrule(l){4-5}\cmidrule(l){6-7}\cmidrule(l){8-9}
			&  CPU &  IT & CPU &  IT  & CPU &  IT & CPU &  IT   \\\specialrule{0em}{1pt}{1pt}\hline\specialrule{0em}{3pt}{3pt}
			
			GMRES & 	2.51E+02 &4696 &	-- &	-- & -- & -- & --	& --
			\\\specialrule{0em}{3pt}{3pt}
			
			CNAS-GMRES & 	1.27E-01 &	16 &	2.66E-01 &	16 & 7.88E-01 & 16 & 1.45E-00  &	18			
			\\\specialrule{0em}{3pt}{3pt}
			
			GE & 2.12E+01 &	N/A  &	1.70E+02  &	N/A  & 1.47E+03  & N/A & 1.98E+04 &  N/A
			\\\specialrule{0em}{3pt}{3pt}	\hline
			
	\end{tabular}}
\end{table}

\begin{table}[htbp]
	\setlength{\abovecaptionskip}{0pt}
	\setlength{\belowcaptionskip}{10pt} \centering{
		\caption{\label{tab:omega alp=1.1} 
			The empirical optimal parameters of CNAS-GMRES  when $\alpha=1.1:0.2:1.9$, $N=200$.}
		\begin{tabular}{llcccc}\specialrule{0em}{2pt}{2pt}\hline\specialrule{0em}{2pt}{2pt}
			& $M$ &  3200 &  6400 &12800  & 25600   \\\specialrule{0em}{1pt}{1pt}\hline\specialrule{0em}{3pt}{3pt}
			
			\multirow{2}{*}{$\alpha=1.1$}  & $\omega_u$  &[0.16,0.24] &[0.16,0.24] &[0.15,0.24] & [0.17,0.24]	
			\\\cmidrule(l){2-6}
			& $\omega_v$ &[0.18,0.25]  &[0.18,0.25] &[0.16,0.22]
			&[0.20,0.25]	
	\\\specialrule{0em}{1pt}{1pt}\hline
            \multirow{2}{*}{$\alpha=1.3$}  & $\omega_u$  &[0.19,0.24] &[0.19,0.24] &[0.18,0.24] & [0.17,0.24]	
			\\\cmidrule(l){2-6}
			& $\omega_v$ &[0.20,0.24]  &[0.20,0.25] &[0.21,0.23]
			&[0.19,0.23]	
			\\\specialrule{0em}{1pt}{1pt}\hline
            \multirow{2}{*}{$\alpha=1.5$}  & $\omega_u$  &[0.09,0.24] &[0.20,0.24] &[0.17,0.24] & [0.17,0.24]	
			\\\cmidrule(l){2-6}
			& $\omega_v$ &[0.10,0.25]  &[0.18,0.25] &[0.17,0.24]
			&[0.18,0.24]	
			\\\specialrule{0em}{1pt}{1pt}\hline
			\multirow{2}{*}{$\alpha=1.7$}  & $\omega_u$  &[0.26,0.34] &[0.26,0.34] &[0.18,0.24] & [0.14,0.24]	
			\\\cmidrule(l){2-6}
			& $\omega_v$ &[0.30,0.43]  &[0.28,0.34] &[0.20,0.25]
			&[0.16,0.25]	
			\\\specialrule{0em}{1pt}{1pt}\hline
			\multirow{2}{*}{$\alpha=1.9$}  & $\omega_u$  &[0.09,0.34] &[0.08,0.34] &[0.06,0.24] & [0.21,0.24]	
			\\\cmidrule(l){2-6}
			& $\omega_v$ &[0.10,0.35]  &[0.09,0.34] &[0.09,0.25]
			&[0.22,0.25]	
			\\\specialrule{0em}{1pt}{1pt}\hline
	\end{tabular}}
\end{table}

\begin{figure}[htbp]
	\centering
	\includegraphics[scale=0.6]{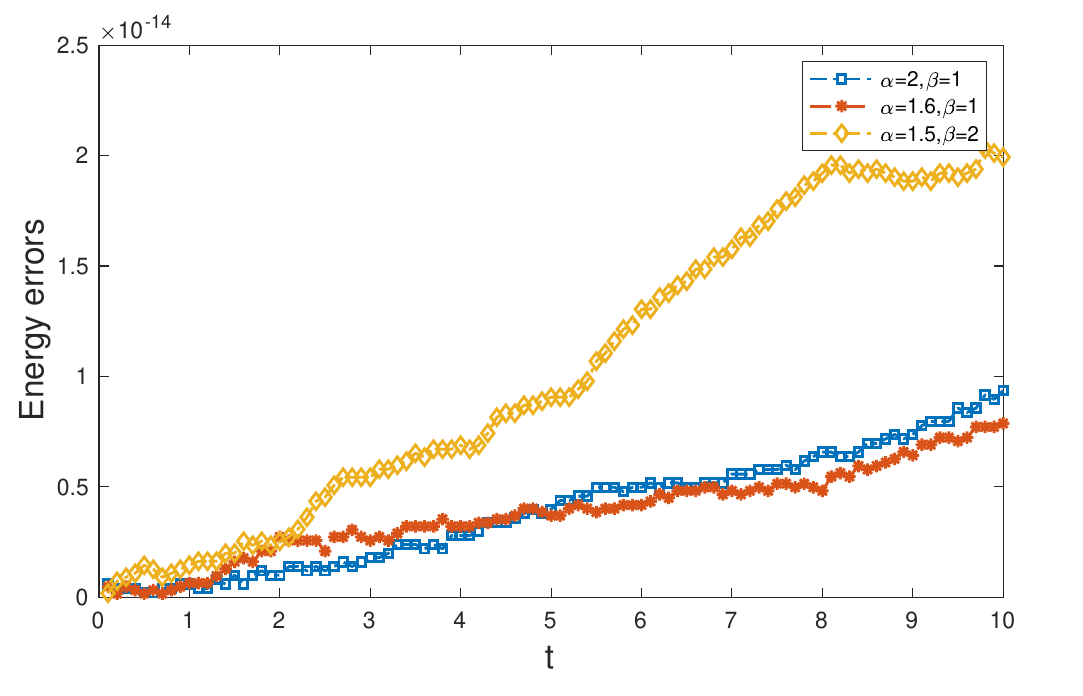}
	\caption{
		The relative errors of the discrete energy, i.e., $|(E^n-E^0)/E^0|$, when $h=0.1$, $\tau = 0.01$.
	}
	\label{fig:energy2}
\end{figure}

Tables \ref{tab:qcon11}-\ref{tab:qcon22} and Figure \ref{fig:energy2} report the relative errors of the discrete mass and  energy. The corresponding space mesh size and time step size are $h=0.1$, $\tau=0.01$. At each time level, the block linear systems related to $u$ and $v$ are solved by CNAS-GMRES, and terminated when the $\ell_2$-norm relative residuals of the block linear systems reduced below $10^{-15}$.   The small errors shown in Tables \ref{tab:qcon11}-\ref{tab:qcon22} and Figure \ref{fig:energy2} show that both the discrete mass and energy of the LICD scheme are conserved by the linear solver CNAS-GMRES in the coupled case with $<\alpha=2,\beta=1>$, $<\alpha=1.6,\beta=1>$, and $<\alpha=1.5,\beta=2>$.

\begin{table}[htbp]
	\renewcommand\arraystretch{1.6}
	
	\tabcolsep=0.25cm
	\setlength{\abovecaptionskip}{2pt}
	\setlength{\belowcaptionskip}{18pt} \centering{
		\caption{
			\label{tab:qcon11}
			The relative errors of the discrete mass in terms of $u$, i.e., $|({Q}_1^{n}-{Q}_1^{0})/{Q}_1^{0}|$, when $h=0.1$, $\tau=0.01$.}
		\begin{tabular}{llllll}
			
			\hline & $t=2$ & $t=4$ & $t=6$ & $t=8$ &$t=10$  \\
			\hline
			
			$\alpha=2,\beta=1$ & $4.4409\mathrm{e}-015$ & $3.8857 \mathrm{e}-015$ & $6.2172\mathrm{e}-015$ & $7.9936\mathrm{e}-015$ &$1.0749\mathrm{e}-014$ \\
			
			$\alpha=1.6,\beta=1$ & $8.8818 \mathrm{e}-016$ & $3.5527\mathrm{e}-015$ & $3.7748 \mathrm{e}-015$ &  $4.2188\mathrm{e}-015$&$3.7748\mathrm{e}-015$ \\
			
			$\alpha=1.5,\beta=2$ & $1.1102\mathrm{e}-015$ & $6.6613 \mathrm{e}-016$ & $2.8866 \mathrm{e}-015$ & $5.5511 \mathrm{e}-015$&$5.3291\mathrm{e}-015$ \\
			\hline
	\end{tabular}}
\end{table}

\begin{table}[htbp]
	\renewcommand\arraystretch{1.6}
	
	\tabcolsep=0.25cm
	\setlength{\abovecaptionskip}{2pt}
	\setlength{\belowcaptionskip}{18pt} \centering{
		\caption{
			\label{tab:qcon22}
			The relative errors of the discrete mass in terms of $v$, i.e., $|({Q}_2^{n}-{Q}_2^{0})/{Q}_2^{0}|$, when $h=0.1$, $\tau=0.01$.}
		\begin{tabular}{llllll}
			
			\hline & $t=2$ & $t=4$ & $t=6$ & $t=8$ &$t=10$  \\
			\hline
			
			$\alpha=2,\beta=1$ & $9.9920\mathrm{e}-016$ & $4.4409 \mathrm{e}-016$ & $2.6645\mathrm{e}-015$ & $3.9968\mathrm{e}-015$ &$9.6589\mathrm{e}-015$ \\
			
			$\alpha=1.6,\beta=1$ & $1.5543 \mathrm{e}-015$ & $4.4409\mathrm{e}-016$ & $3.2204 \mathrm{e}-015$ &  $1.3323 \mathrm{e}-015$&$1.3323\mathrm{e}-015$ \\
			
			$\alpha=1.5,\beta=2$ & $4.4409 \mathrm{e}-016$ & $2.2204 \mathrm{e}-016$ & $1.1102 \mathrm{e}-015$ & $4.2188 \mathrm{e}-015$&$3.5527\mathrm{e}-015$ \\
			\hline
	\end{tabular}}
\end{table}

\begin{figure}[htbp]
	\centering
	\subfloat{\includegraphics[scale=0.38]{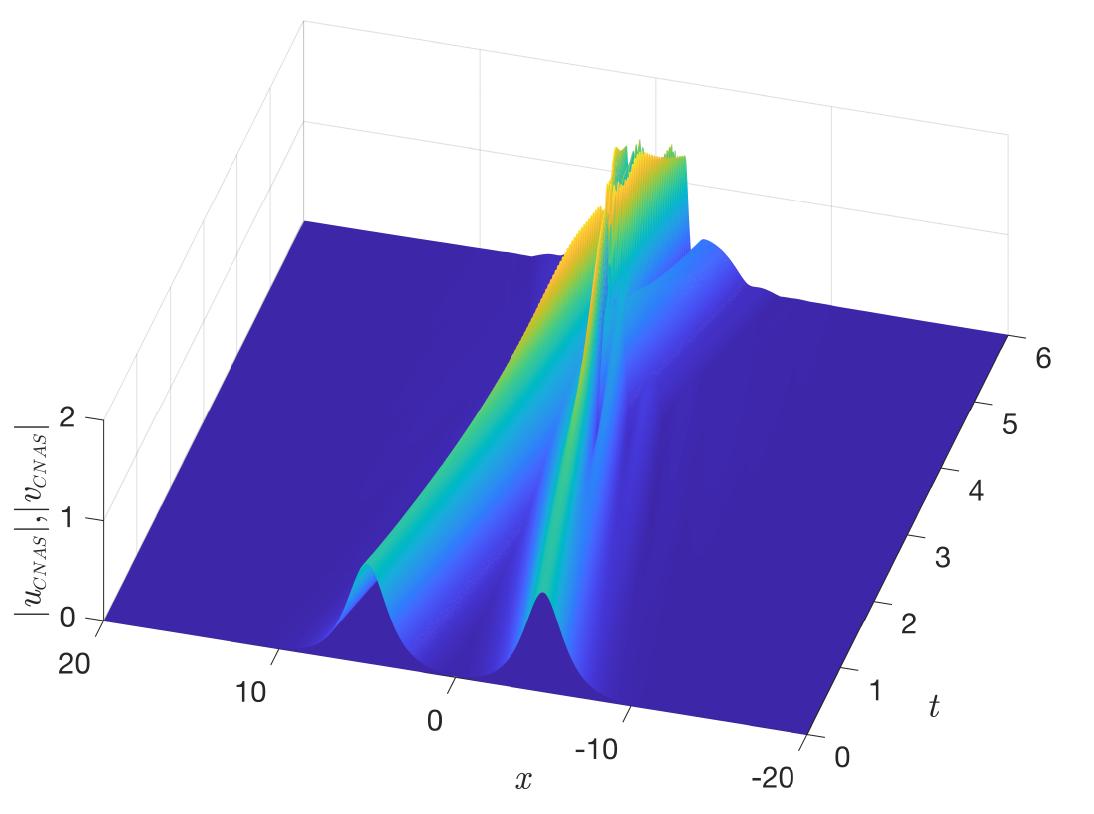}}
	\subfloat{\includegraphics[scale=0.38]{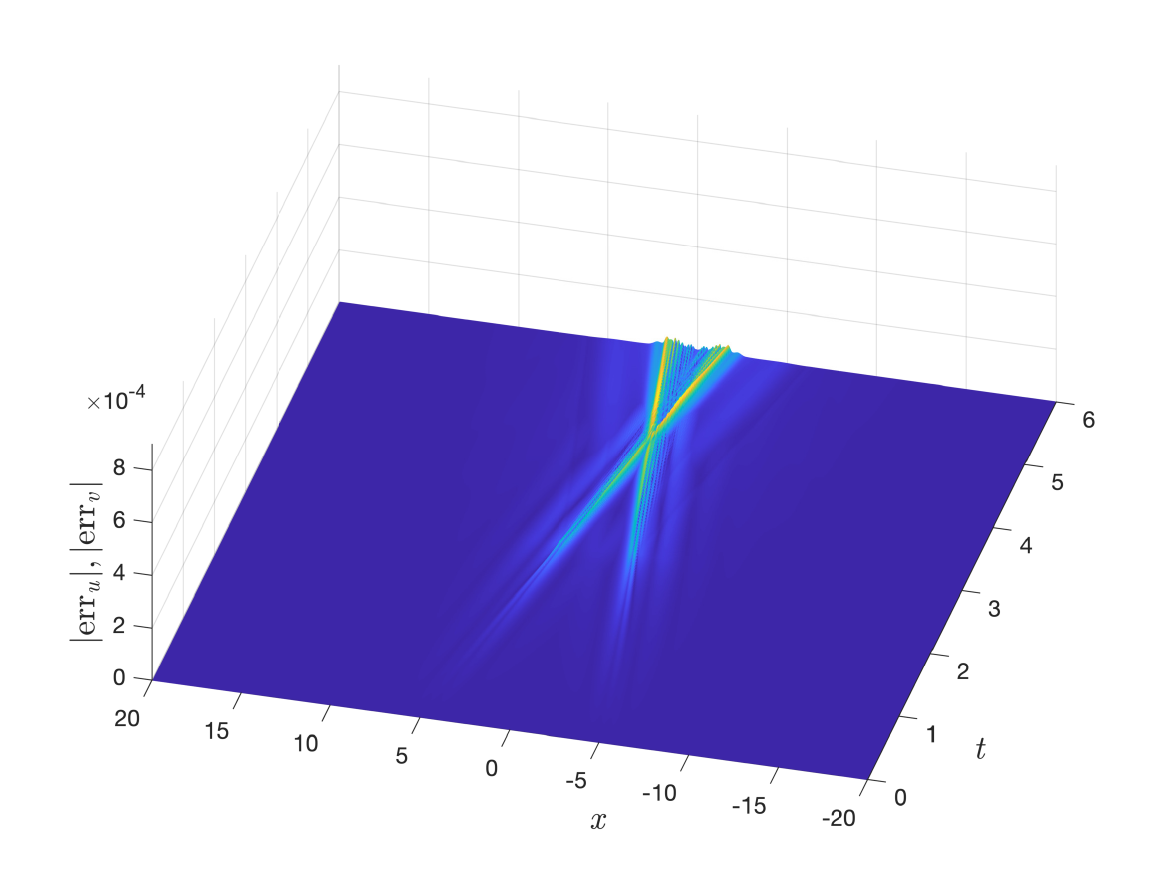}}
	\caption{
		The numerical solutions (left) and their errors (right)  of  the space fractional CNLS equations (\ref{equ:CNLS}) when $\alpha=1.1$, $M=800$, $N=600$.}
	\label{figCsingal_1.1}
\end{figure}

\begin{figure}[htbp]
	\centering
	\subfloat{\includegraphics[scale=0.38]{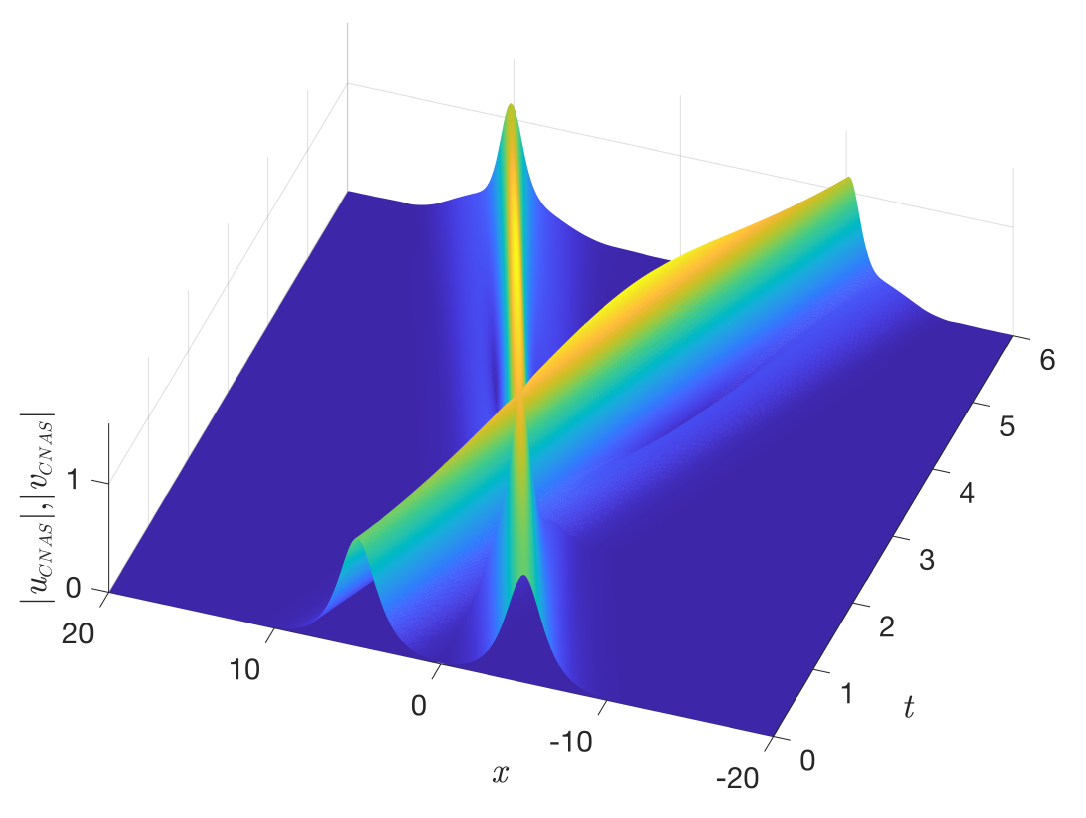}}
	\subfloat{\includegraphics[scale=0.38]{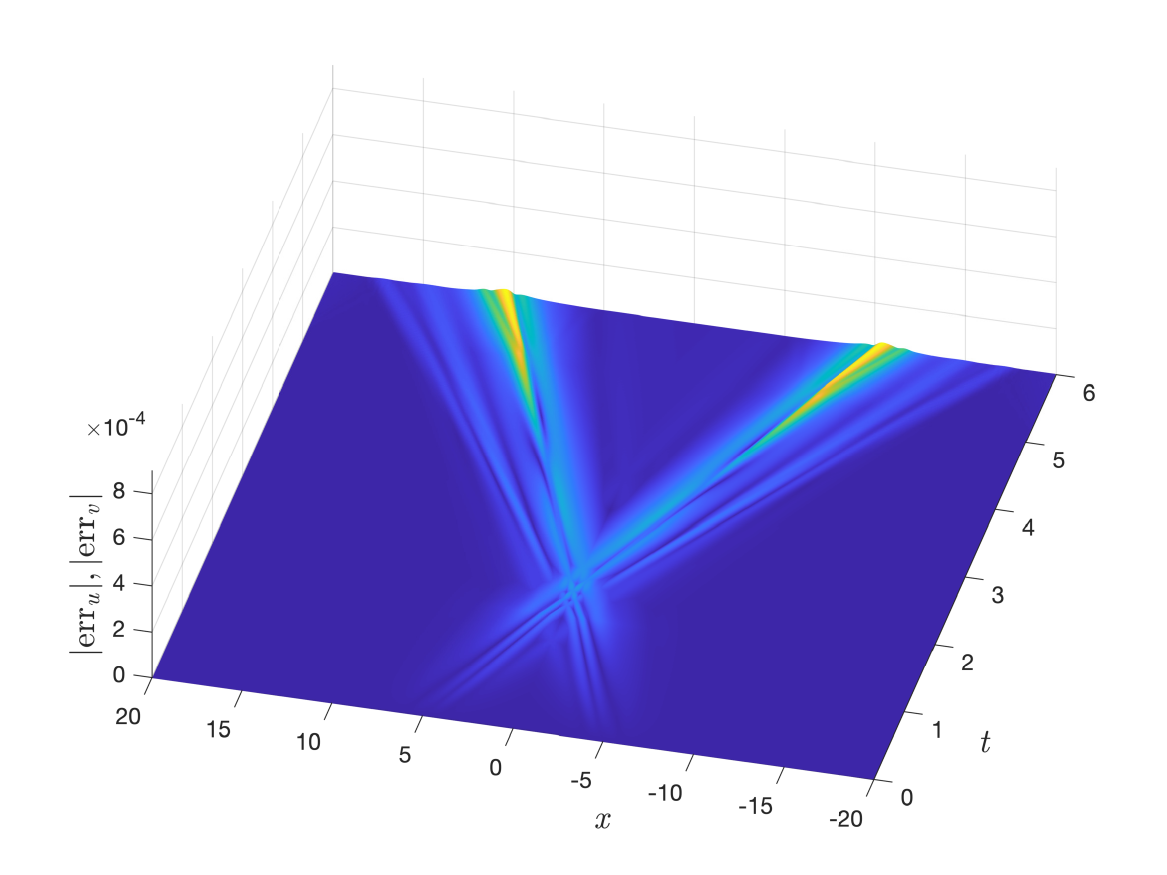}}
	\caption{
		The numerical solutions (left) and their errors (right)  of  the space fractional CNLS equations (\ref{equ:CNLS}) when $\alpha=1.5$, $M=800$, $N=600$.}
	\label{figCsingal_1.5}
\end{figure}

\begin{figure}[htbp]
	\centering
	\subfloat{\includegraphics[scale=0.38]{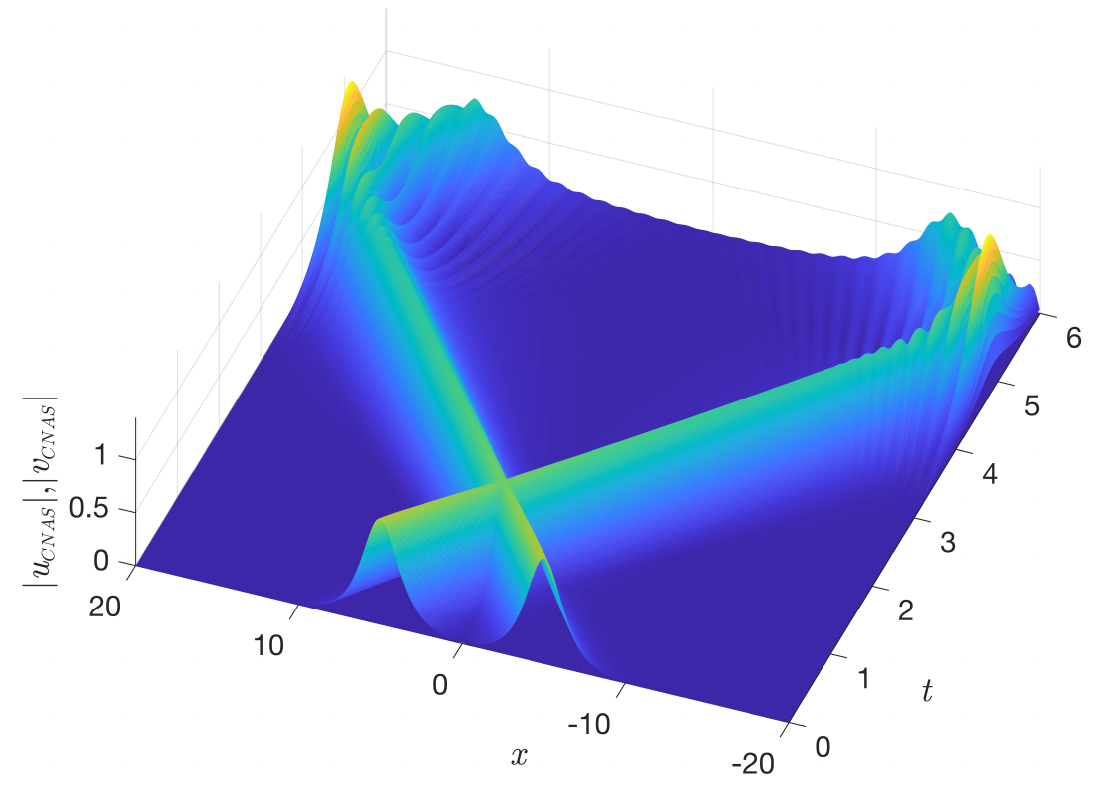}}
	\subfloat{\includegraphics[scale=0.38]{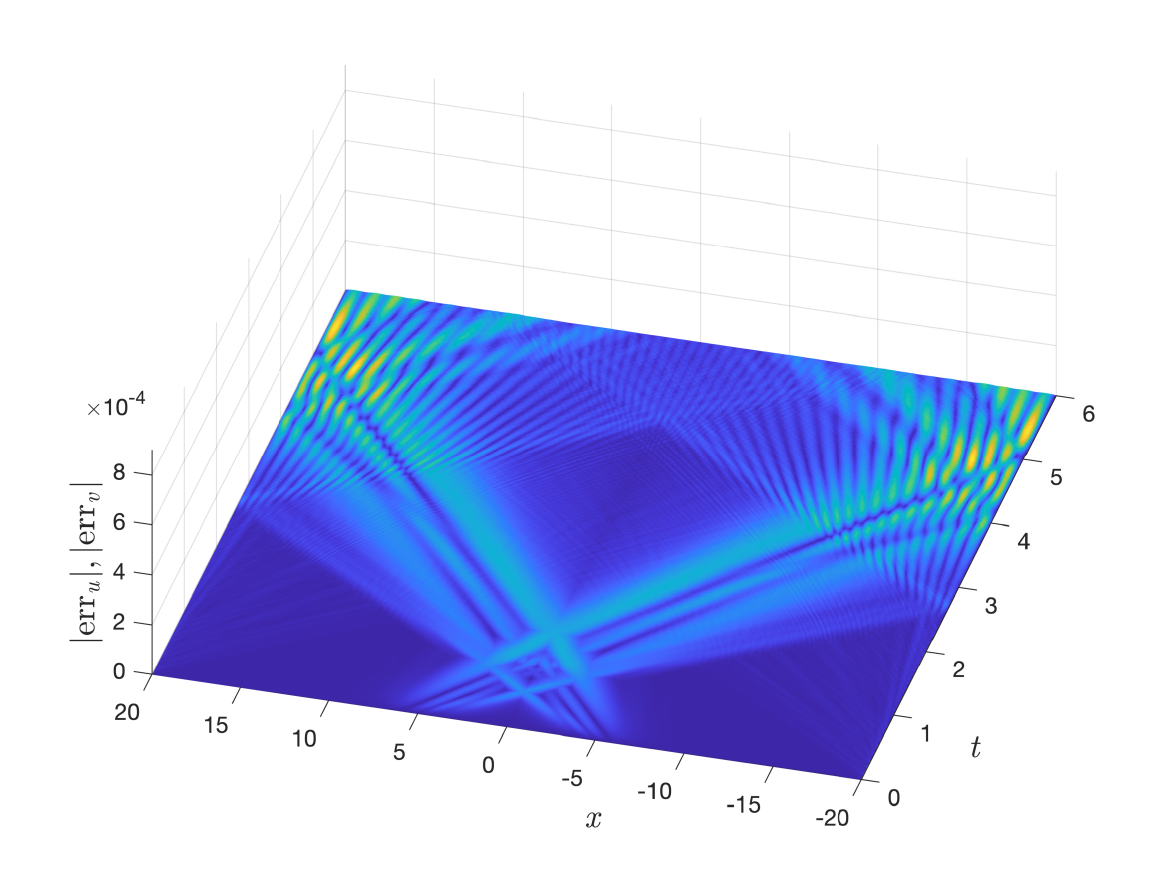}}
	\caption{
		The numerical solutions (left) and their errors (right)  of the space fractional CNLS equations (\ref{equ:CNLS}) when $\alpha=1.9$, $M=800$, $N=600$.}
	\label{figCsingal_1.9}
\end{figure}

\begin{figure}[htbp]
	\centering
	\subfloat{\includegraphics[scale=0.38]{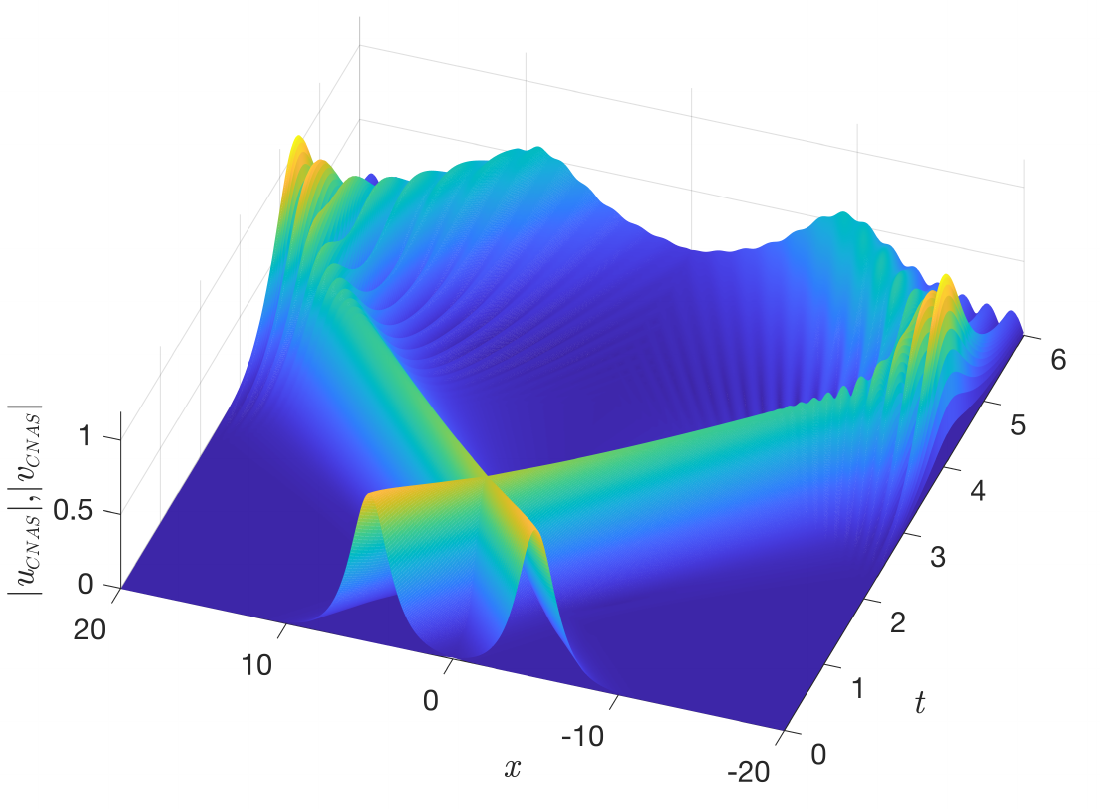}}
	\subfloat{\includegraphics[scale=0.38]{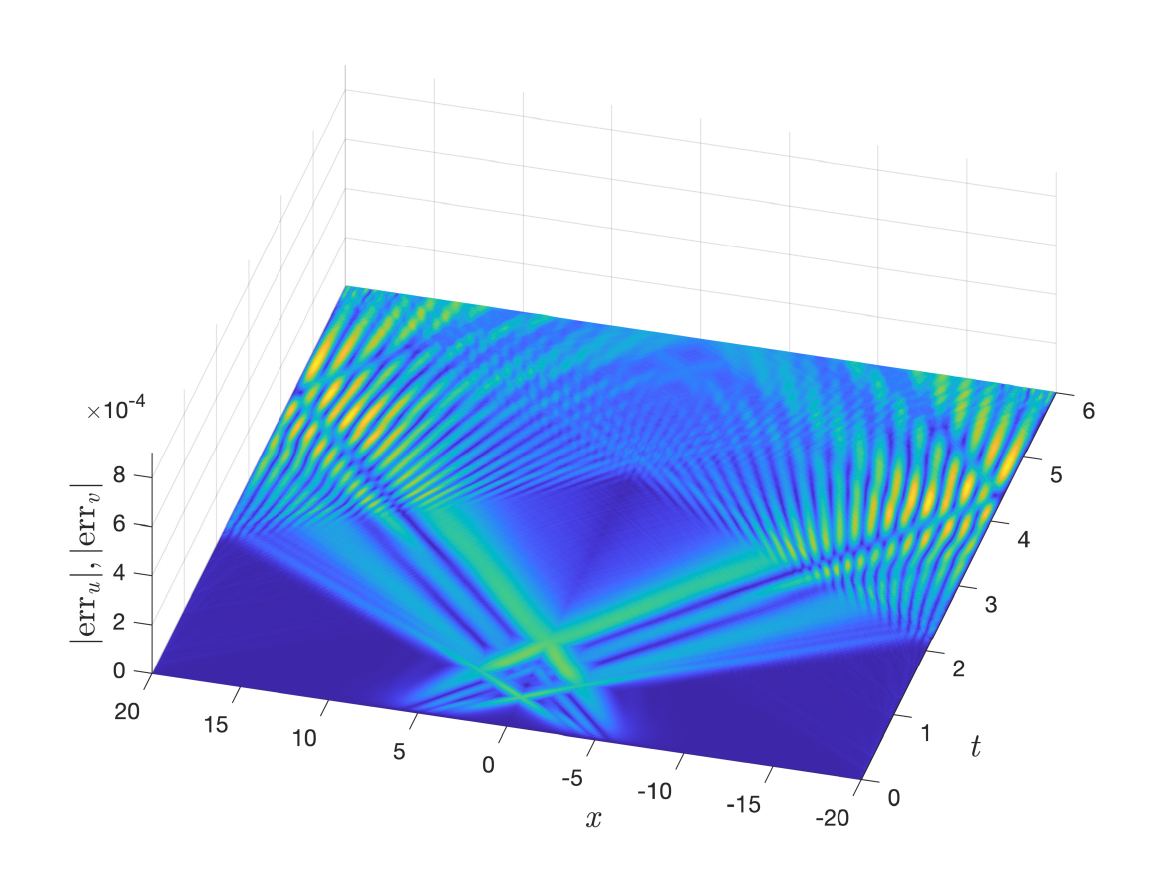}}
	\caption{
		The numerical solutions (left) and their errors (right)  of  the space fractional CNLS equations (\ref{equ:CNLS}) when $\alpha=2$, $M=800$, $N=600$.}
	\label{figCsingal_2}
\end{figure}

Figures \ref{figCsingal_1.1}-\ref{figCsingal_2} depict the numerical solutions (left), i.e., $u_{\text{\tiny CNAS}}$ and $v_{\text{\tiny CNAS}}$, obtained by  CNAS-GMRES, and their errors (right), i.e., $\text{err}_{u} = |u_{\text{\tiny CNAS}}-u_{\text{\tiny GE}}|$ and $\text{err}_{v} =|v_{\text{\tiny CNAS}}-v_{\text{\tiny GE}}|$
with the exact solutions $u_{\text{\tiny GE}}$ and $v_{\text{\tiny GE}}$ of the LICD scheme (i.e., the solutions $u_{\text{\tiny GE}}$ and $v_{\text{\tiny GE}}$ obtained by GE) for the space fractional CNLS equations (\ref{equ:CNLS}) when $\alpha=1.1:0.4:1.9$ and $\alpha=2$, $M=800$, $N=600$. The shape of  the wave fronts varies with the fractional order $\alpha$, and it converges to the wave fronts of the standard CNLS equations as $\alpha$ approaches to $2$. The fractional order $\alpha$ affects the collision time of the wave fronts. The larger the value of $\alpha$, the faster the wave fronts move, and the earlier the collision occurs. It can be seen from Figures \ref{figCsingal_1.9}-\ref{figCsingal_2} that reflections occur after the wave fronts reach the boundary of the space-time domain. Obviously, no wave front reflection will occur when the space interval is not truncated. In addition, the errors of the numerical solutions stay very small meaning that the linear solver CNAS-GMRES is reliable.

\section{The concluding remarks}\label{conclu}
We focus on designing effective and efficient linear solvers for complex symmetric and indefinite linear systems of the form $(D-T+\imath I) \textbf{u} = \textbf{b}$, which have Toeplitz-plus-diagonal and complex symmetric structure. These linear systems are originated from 1D space fractional Schrödinger equations with attractive interaction of particles discretized by the LICD scheme. Specifically, the NASS iteration method and the naturally induced NASS preconditioner are proposed. An easy to implement and extremely efficient preconditioner is constructed by accurately approximating the involved Toeplitz matrix block by a circulant matrix, which is called the CNAS preconditioner. Theoretically, the new iteration method and preconditioners are analyzed in detail. The CNAS preconditioned GMRES method is verified to be an effective and efficient linear solver for $(D-T+\imath I) \textbf{u} = \textbf{b}$ by numerical experiments based on 1D space fractional CNLS equations. However, the construction of the new linear solver is due to the constant coefficients of the 1D problem and the adopted uniform grids. In the future, we can extend the NASS iteration method and the related preconditioner to higher dimensional problems. Moreover, when variable coefficients appear in the problem, no explicit Toeplitz-plus-diagonal structure can be found in the derived linear systems, finding possible implicit structure and extending the techniques in this paper to construct efficient solvers may be a great challenge. Finally, when the problem is discretized on non-uniform grids, the Toeplitz-plus-diagonal structure of the derived linear system may be completely lost, a possible way to construct fast solvers is to combine the hierarchical-matrix approach \cite{BM2008SpringerBook,HKL2016CPAM} and the framework proposed in this paper.

\section*{Acknowledgments}
This work was funded by the National Natural Science Foundation (No. 11101213 and No. 12071215), China.

\end{document}